\documentclass{amsart}
\usepackage[utf8]{inputenc}

\usepackage[english]{babel}

\usepackage{amssymb}
\usepackage{amsmath}
\usepackage{amsmath,amscd}
\usepackage{amsthm}
\usepackage{enumerate}
\usepackage{cite}
\usepackage{xcolor}
\usepackage[final]{listings}
\usepackage[section]{algorithm}

\makeindex

\title[On the description of identifiable quartics]
{On the description of identifiable quartics}
\date{}

\newcommand{\C}{\mathbb{C}}
\newcommand{\Z}{\mathbb{Z}}
\newcommand{\Pj}{\mathbb{P}}

\newcommand{\N}{\mathbb{N}}

\newcommand{\Oc}{\mathcal{O}}
\newcommand{\Ic}{\mathcal{I}}
\newcommand{\Ec}{\mathcal{E}}
\newcommand{\Zc}{\mathcal{Z}}
\newcommand{\Tc}{\mathcal{T}}
\newcommand{\Bc}{\mathcal{B}}

\newtheorem{thm0}{Theorem}[section]
\newtheorem{prop0}[thm0]{Proposition}

\newtheorem{coro0}[thm0]{Corollary}

\theoremstyle{definition}
\newtheorem{defn0}[thm0]{Definition}

\newtheorem{exa0}[thm0]{Example}

\newtheorem{rem0}[thm0]{Remark}


\author[E.~Angelini]{Elena Angelini}
\address{Elena Angelini: Dipartimento di Ingegneria dell'Informazione e Scienze Matematiche, Universit\`a di Siena, Italy}
\email{elena.angelini@unisi.it}

\author[L.~Chiantini]{Luca Chiantini}
\address{Luca Chiantini: Dipartimento di Ingegneria dell'Informazione e Scienze Matematiche, Universit\`a di Siena, Italy}
\email{luca.chiantini@unisi.it}
\thanks{The authors are members of the Italian GNSAGA-INDAM}

\begin{document}

\begin{abstract}
In this paper we study the identifiability of specific forms (symmetric tensors), with the target of extending recent methods for the case of 
$3$ variables to more general cases.
In particular, we focus on forms of degree $4$ in $5$ variables. By means of tools coming from classical algebraic geometry, 
such as Hilbert function, liaison procedure and Serre's construction, we give a complete geometric description and criteria of identifiability 
for ranks $\geq 9$, filling the gap between rank $\leq 8$, covered by Kruskal's criterion,
and $15$, the rank of a general quartic in $5$ variables. For the case $r=12$, we construct an effective algorithm that guarantees that a given decomposition
is unique.
\end{abstract}

\maketitle

\section{Introduction}

We consider the study of symmetric tensors $T$ over $\C$, which we will identify with homogeneous polynomials (forms),
with respect to their Waring (i.e. symmetric) rank and identifiability.
Our point of view is the following. Assume that a specific form $T$ is given and assume we know an {\it expression} of $T$ as a sum of $r$ powers of linear forms,
\begin{equation}\label{expr} T = \lambda_1T_1+\dots+\lambda_rT_r. \end{equation}
Then:

({\bf Q1}) is the length $r$ minimal? In other words, is $r$ the {\it rank} of $T$?\\
If the answer to the previous question is positive, then another question arises: 

({\bf Q2}) is \eqref{expr} the {\it unique} expression of minimal length for $T$?\\
Here {\it unique} means, of course, up to a permutation of the summands. When uniqueness holds, we will say that $T$ is {\it identifiable}.

 The most celebrated criterion that answers the questions is due to Kruskal (see \cite{Kruskal77}). It works provided
 that an inequality involving some invariants of the expression (the {\it Kruskal's ranks}) is satisfied. Kruskal's criterion, which indeed works
 for general multilinear tensors, is continuously employed in applications of tensor analysis to statistics, signal theory, chemistry,
 quantum information theory, artificial intelligence, etc. We mention, for one, the paper \cite{RaoLiZhang18}.  
  It is known that the Kruskal's criterion is sharp. Nothing better can be done if one uses solely the Kruskal's ranks of expression \eqref{expr} (see \cite{Derksen13}). 
Since Kruskal's inequality cannot hold outside a precise range for the length $r$, applications of the Kruskal's criterion are bounded
 to some, tipically quite small, values of $r$. The range in which Kruskal's criterion possibly applies is much smaller, in fact,
 than the  range in which identifiability holds for general forms (see \cite{COttVan17b}).
 
Actually, there are methods that determine the identifiability of $T$ in a wider range for $r$, provided that
 one computes some higher invariants of  \eqref{expr}. For instance, the flattening procedure 
 can determine the identifiability of $T$  even when Kruskal's inequality fails (see \cite{Domanov, DomaLath15}). 
 Yet, it seems hard to cover with flattening methods the whole range in which identifiability can hold.
 
  In a series of papers \cite{BallBern12a, Ball19, COttVan17b, AngeCVan18, AngeCMazzon19, AngeC20, AngeC}, it is shown that a geometric approach can
determine the identifiability of tensors in the whole range in which the property can hold. 
  The new  method starts by considering a finite set $A=\{P_1,\dots,P_r\}$ of points in a projective space ($\Pj^n$, for forms in $n+1$ variables),
 naturally associated to   \eqref{expr}. The analysis of the geometry of $A$
 can exclude the existence of a second expression for $T$ of length $\leq r$, for $r$ ranging in a wide set of values. 
 
In the new methods, properties of the finite set $A\subset \Pj^n$ needed to examine the identifiability of $T$
 are rather deep and require advanced tools, both from algebra and geometry. They are based on
 an accurate analysis of the   the resolution of the homogeneous ideal, the Hilbert function and the Hilbert-Burch matrix of $A$,
 together with liaison techniques. For ternary forms, when the set $A$ lives in $\Pj^2$, these properties 
describe the situation quite completely. When the number of variables grows, and $A$ lives in higher 
 dimensional projective spaces, even deeper theoretical results are necessary to understand the identifiability of specific forms.
 Not all the geometric tools necessary for the analysis are presently available. This turns out to be a stimulating
 challenge for algebraic geometry, that could produce substantial advances both in the theory and for applications. 
 
 We produce a pattern that introduces some deeper geometric tools (minimal resolution conjecture, vector bundles on surfaces)
and determines a criterion for the identifiability of some forms $T$ in more than three variables. These advanced  tools
are not new, even if rather  recent. What is new is their combination in a procedure that excludes the existence of sets of projective points forming an 
alternative decomposition of $T$.  We focus in particular on the case of quartics in $5$ variables, even if a similar theory can handle  more general forms. 
The identifiability of quartics in $2,3,4$ variables turns out to be widely understood (for $4$ variables, see e.g. \cite{AngeCVan18}).
The case of quartics in $5$ variables was not completely covered in the previous literature.

For quartics in $5$ variables we will study in details the uniqueness of expression \eqref{expr} when $r$ is equal to $9,10,11,12,13$.
 Notice that the range $r<9$ is covered by the original Kruskal's criterion.
For higher values of $r$, notice that $r=14$ is a special case of the Alexander-Hirschowitz theorem \cite{AlexHir95},
while $15$ is the generic value for the rank. It follows that for $r\geq 14$
the standard map from the abstract secant variety to the secant variety of the $4$-Veronese variety of $\Pj^4$ (as in Remark \ref{stacked} below)
 has positive dimensional fibers, and
this implies (see e.g. \cite{CCi02a}) that an expression of length $r\geq 14$ cannot be 
 unique. Indeed, if one such expression exists, then infinitely many expressions of the same length must exist.
So, identifiability is excluded for $r\geq 14$.
Thus we have a  method which  tests the identifiability of $T$ for all possible values of  $r$.
 
Case  $r\leq 11$ turns out to be different from cases $r=12,13$. For $r\leq 11$, when  $A$ is sufficiently general 
 (in a  precise sense, that can be tested by computer algebra  algorithms),
 then identifiability holds, and no forms spanned by  the powers $T_1,\dots,T_r$ of \eqref{expr}, except trivially those
 spanned by proper subsets of the $T_i$'s, can have alternative expressions of length $\leq r$. 
 
Conversely for $r=12$, even if the $T_i$'s are
 general, yet their linear span $L$ contains some special tensors with a second, different decomposition of the same length. So, one needs a sharp
 criterion to exclude the existence of a second expression for a given $T$. Our method distinguishes between different forms in $L$, 
 i.e. it takes care of the (ratio between)  coefficients  in the expressions  $a_1T_1+\dots+a_rT_r$ of forms lying in  $L$.
 Following our method, we construct an algorithm which tests if a given decomposition of length $12$ of a specific form $T$
 is unique or not (see Section 3.3). The algorithm is {\it effective} in the sense of \cite{COttVan17b}: it will  return a positive answers for all $T$
 outside a (Zariski closed) subset of measure $0$ in the space of forms. An implementation is described in the file {\tt ancillary.txt}.

The case $r=13$ is even more involved. Not only the span $L$ of general $T_i$'s  contains forms with alternative expressions of length
$13$, but it also contains forms of rank $12$ (not generated by proper subsets of  $\{T_1,\dots,T_{13}\}$). Even in this case,
we are able to produce a procedure that excludes the existence of alternative expressions of length $\leq 13$ for a specific $T$.

In conclusion, we show that advanced geometric tools give a method to test the uniqueness (and minimality) of an
expression  \eqref{expr} for quartics in $5$ variables. The criterion is effective, in the sense of \cite{COttVan17b}: it will
give a (positive) answer for forms which lie outside an algebraic (Zariski closed) subset in the space of forms $S^4(\C^5)$.

The structure of the paper is the following. Section \ref{sec:notation} is devoted to the main notation, concepts and results, used throughout the paper, 
coming both from tensors setting (such as Kruskal's criterion and its generalizations) and from classical algebraic geometry 
(Hilbert function, Cayley-Bacharach property, Liason procedure and Serre's construction). 
By means of the above mentioned fundamental tools, in Section \ref{sec:quartics} we give a  geometric 
procedure to test the identifiability for quartics in five variables, complete in a dense subset, which covers also the cases of rank 
$11,12,13$ that, at the best of our knowledge,  were not  covered in the mathematical literature on this subject.

\section{Preliminaries}\label{sec:notation}

\subsection{Notation}
For $d,n \in \N$, let $ \C^{n+1} $ be the space of linear forms in $ x_{0}, \ldots, x_{n} $, thus $ S^{d} \C^{n+1}$ 
is the space of forms of degree $d$ in $n+1$ variables over $\C$. Every $ T \in S^{d} \C^{n+1} $ defines an element of 
$ \Pj(S^{d} \C^{n+1}) \cong \Pj^{N} $ $( N = \binom{n+d}{d} - 1) $, which we still denote by $T$. Moreover $ v_{d}: \Pj^{n} \rightarrow \Pj^{N} $ 
is the \emph{Veronese embedding} of $ \Pj^{n} $ of degree $ d $, i.e. 
$$ v_{d}([a_{0}x_{0}+ \ldots + a_{n}x_{n}]) = [(a_{0}x_{0}+ \ldots + a_{n}x_{n})^{d}]. $$ 
For any finite set $ A = \{P_{1}, \ldots, P_r\} \subset \Pj^{n} $, $ \langle v_{d}(A) \rangle $ is the linear space in $\Pj^N$ spanned by the points 
$ v_{d}(P_{1}), \ldots, v_{d}(P_{r}) $. The cardinality of $ A $ is usually denoted by $ \ell(A) $ and $ I_{A} $, $ \mathcal{I}_{A} $ are, respectively, 
the ideal of $ A $ in the polynomial ring $ R = \C[x_{0}, \ldots, x_{n}] $ and the ideal sheaf of $ A $ on $ \Pj^{n} $. 
We will denote with subscripts the homogeneous pieces of $R$ and its ideals.

\begin{defn0}
Given a finite set $ A \subset \Pj^n $ and a form $ T \in S^{d} \C^{n+1} $, we say that:
\begin{itemize}
\item[$\bullet$] $ A $ \emph{computes} $ T $, or that $A$ is a decomposition of $T$, if $ T\in \langle v_{d}(A) \rangle$;
\item[$\bullet$] $ A $ is \emph{non-redundant} if $ A $ computes $ T $ and there are no proper subsets $ A' $ of $ A $ 
such that $A'$ computes $T$;
\item[$\bullet$] $ A $ is \emph{minimal} if $ A $ computes $ T $ and there are no sets $B$, with $\ell(B)<\ell(A)$, 
such that $B$ computes $T$;
\item[$\bullet$] $\ell(A)$ is the (Waring) \emph{rank} of $T$ if $ A $ is minimal;
\item[$\bullet$] $T$ is \emph{identifiable} if $ A $ is the unique set such that $ \ell(A) $ equals the rank of $T$. 
\end{itemize} 
\end{defn0}

\subsection{Kruskal's criterion for forms and its extensions}\label{sec:Kr}

In the mathematical literature, one of the most famous criteria for detecting the identifiability of a tensor is due to Kruskal. It is based on the concept of $ d $-th Kruskal's rank, $ k_{d} $, of a finite set, for whose definition and main properties we refer to $\S$ 2.2 of \cite{AngeC}. Here we recall one extension of this criterion, adapted to the case of forms.

\begin{thm0}[Reshaped Kruskal's Criterion, see \cite{COttVan17b}]\label{thm:kr}
Assume $ d \geq 3 $ and let $ A \subset \mathbb{P}^{n} $ be a non-redundant decomposition of $ T \in \Pj(S^{d} \C^{n+1}) $.  
Fix a partition  $ d = d_{1}+d_{2}+d_{3} $ with $ d_{1} \geq d_{2} \geq d_{3} \geq 1 $ and denote by $ k_{d_{i}}(A) $ the $ d_{i} $-th Kruskal's rank of $ A $. If 
\begin{equation}\label{eq:Kr}
\ell(A) \leq \frac{k_{d_{1}}(A)+k_{d_{2}}(A)+k_{d_{3}}(A)-2}{2}
\end{equation}
then $ T $ has rank $ \ell(A) $ and it is identifiable.
\end{thm0}

In the case of ternary forms, Theorem \ref{thm:kr} has been recently extended in \cite{AngeC}, \cite{AngeC20}. \\
For the case of quartics in $ n+1 $ variables, we refer to the following extension (see Section 6 of \cite{AngeCVan18}):

\begin{thm0}\label{thm:krquartics} Let $ T \in \Pj(S^{4} \C^{n+1}) $ and let $ A = \{P_{1}, \ldots, P_{2n+1}\} \subset \mathbb{P}^{n} $ be a non-redundant 
decomposition of $ T $. Write $X=v_4(\Pj^n)$ for the image of the Veronese embedding. If
\begin{itemize}
\item[a)]  $ \dim \langle v_{4} (P_{1}), \ldots, v_{4}(P_{2n+1}) \rangle = 2n+1 $,
\item[b)] $ k_{1} (A) = n+1 $,
\item[c)] the linear span of the union of tangent spaces $\bigcup\bold{T}_{X,v_{4} (P_i)}$ has the (expected) dimension $2n^2 + 3n + 1 $,
\end{itemize}
then $T$  is identifiable of rank $ 2n+1 $.\end{thm0}

\begin{rem0}\label{stacked}
Condition c) in the previous statement is linked to the Terracini's Lemma which describes the tangent space to secant varieties. For specific decompositions
$A$ of $T$, the condition can be verified by a linear algebra algorithm introduced in  \cite{COttVan14}.

It is a standard consequence of the Terracini's construction that if $T$ has infinitely many decompositions, then condition c) cannot hold (see e.g.
Lemma 6.5 of \cite{COttVan14} or Proposition 6 of  \cite{Mazzon20}). 
\end{rem0}

\subsection{Hilbert function, first difference, h-vector of finite sets}\label{sec:hilb}

Other important tools for our analysis, which come from classical Algebraic Geometry, are the Hilbert function, its first difference and, consequently, the $ h$-vector of a finite set. For completeness, we briefly recall their definitions.

\begin{defn0}
Let $ Y $ be a set of homogeneous coordinates for a finite set $ Z \subset \mathbb{P}^{n} $.
\begin{itemize}
\item[$\bullet$] The \emph{Hilbert function} of $ Z $ is the map
$$ h_{Z}: \mathbb{Z} \longrightarrow \mathbb{N} $$
such that 
\begin{equation*}
h_{Z}(j) =
\begin{cases} \begin{array}{ll}
0 &  \mbox{for j $<$ 0}\\
rank (ev_{Y}(j)) & \mbox{for j $\geq$ 0}
\end{array}\end{cases}
\end{equation*}
where $ ev_{Y}(j): S^{j}\mathbb{C}^{n+1} \longrightarrow \mathbb{C}^{\ell(Z)} $ is the \emph{evaluation map} of degree $ j $ on $ Y = \{Y_{1}, \ldots, Y_{\ell(Z)}\} $, i.e. the linear map given by $ ev_{Y}(j)(F) = (F(Y_{1}), \ldots, F(Y_{\ell(Z)})). $
\item[$\bullet$]  The \emph{first difference of the Hilbert function} $ Dh_{Z} $ is given by
$$ Dh_{Z}(j) = h_{Z}(j)-h_{Z}(j-1), \,j \in \Z . $$
\item[$\bullet$] The $ h $-vector of $ Z $ is the vector consisting of all the non-zero values of $ Dh_{Z} $.
\end{itemize}
\end{defn0}

We refer to $ \S $ 2.3 of \cite{AngeC} for a list of many useful properties.\\
We just point out  that when $h_Z(d)=\ell(Z)$ we say that $Z$ is separated in degree $d$. Notice that $Z$ is separated in degree $d$ when the evaluation map surjects.
This implies that $Z$ is also separated in any degree $\geq d$.
We stress, for reference, the following standard fact.

\begin{prop0} \label{restri} Assume that $Z$ is separated in degree $d-1$. Then for a general linear linear form $\Lambda\in R$ the restriction map
$(I_Z)_d \to (R/\Lambda)_d$ surjects.
\end{prop0}
\begin{proof} We can identify $R/\Lambda$ with the polynomial ring  in $n$ variables. The restriction map $R_d\to (R/\Lambda)_d$ surjects. The conclusion follows from the snake lemma applied to the diagram
$$\begin{matrix}
0 & \to & (I_Z)_{d-1} &\to  &  R_{d-1} & \to & \mathbb{C}^{\ell(Z)} &\to& 0   \\
 &&   \downarrow \Lambda & &  \downarrow \Lambda & & ||  & &     \\
0 & \to & (I_Z)_d &\to  &  R_d & \to & \mathbb{C}^{\ell(Z)} &\to& 0\\
 &&   & &  \downarrow  & &  & &      \\
& &   & &  (R/\Lambda)_d  & & & & 
\end{matrix}$$
where the leftmost vertical maps are multiplication by $\Lambda$.
\end{proof}

The applications of the Hilbert function to the identifiability of forms is based on the following well-known proposition (see e.g. Lemma 1 of \cite{BallBern12a}).

\begin{prop0}\label{d+1}  
Let $ T \in S^{d}\C^{n+1} $ and let $ A, B \subset \Pj^{n} $ be non-redundant decomposition of $ T $. Then $ Dh_{A \cup B}(d+1) > 0 $.
\end{prop0}

\subsection{The CB property for finite sets}\label{sec:CB} The Cayley-Bacharach property is fundamental to detect all possible $h$-vectors of $ A \cup B $, 
where $ A $ and $ B $ are disjoint non-redundant decompositions of  a form $ T $. 

\begin{defn0}\label{def:CB}
A finite set $Z\subset \Pj^n$ satisfies the \emph{Cayley-Bacharach property in degree $i$}, for simplicity $\mathit{CB}(i)$, 
if $ H^{0}(\mathcal{I}_{Z\setminus\{ P\}}(i)) = H^{0}(\mathcal{I}_{Z}(i)) $ for all $P \in Z$.
\end{defn0}

The application of the $CB$ property to our analysis is based on the following results. For the proofs we refer, respectively, to \cite{AngeCVan18} and \cite{AngeC20}.

\begin{thm0}\label{GKRext}
Let $ Z \subset \Pj^{n} $ be a finite set satisfying $\mathit{CB}(i)$ and let $ j \in \{0, \ldots, i+1 \} $. Then
$$ Dh_{Z}(0)+Dh_{Z}(1)+\cdots + Dh_{Z}(j) \leq Dh_{Z}(i+1-j)+\cdots +Dh_{Z}(i+1).$$
\end{thm0}

\begin{thm0}\label{CBconseq}
Let $ T \in S^{d}\C^{n+1} $ and let $ A, B \subset \Pj^{n} $ be non-redundant decompositions of  $ T $ such that $ A \cap B = \emptyset $. 
Then $ A\cup B $ satisfies $\mathit{CB}(i)$, for any $ i \in \{0, \ldots, d\} $.
\end{thm0}

\subsection{Liaison and mapping cone}\label{sec:liaison} Several constructions will be based on the notion of {\it liaison}, or linkage,  of finite sets,
that we recall here briefly. We point to \cite{Ferrand75}, \cite{PeskineSzpiro74}, and \cite{Migliore} for details and proofs.

\begin{defn0} We say that two finite sets $A,B\subset\Pj^n$ are \emph{linked} when there exists a complete intersection $Z$ such that
$I_B=I_Z:I_A$. When $A,B$ are disjoint, this simply means that $A\cup B=Z$. When $A,B$ are linked by $Z$, then we also say that $B$ is the \emph{residue} of $A$
with respect to $Z$.
\end{defn0}

For a finite set  $ Z $ a \emph{ resolution} of the ideal $I_Z$ is an exact sequence of free modules and degree $0$ maps:
$$ 0\to F_n \to F_{n-1} \to \dots \to F_1\to I_Z\to 0.$$
When $Z$ is a complete intersection, a resolution is given by the Koszul complex (see \cite{Eisenbud}, section 17), so that $F_n$ has rank $1$.\\
If $A\subset Z$, then the resolutions of the ideals of $A,Z$ determine a commutative diagram of resolutions:
$$
\begin{matrix}
0 &\to& F_n &\to& F_{n-1}&\to& \dots &\to & F_1&\to&  I_Z  &\to& 0 \\ 
& & \vspace{.1cm} & & & & & & & &  & &  \\

 & & \big\downarrow & & \big\downarrow & & & & \big\downarrow & & \big\downarrow  & &  \\

& & \vspace{.1cm} & & & & & & & &  & &  \\

0 &\to& G_n &\to& G_{n-1}&\to& \dots &\to & G_1&\to&  I_A &\to& 0
\end{matrix}$$
where the rightmost vertical arrow is the natural inclusion.\\
If furthermore $Z$ is a complete intersection of type $d_1,\dots,d_n$, then a (non-necessarily minimal) resolution of the ideal of the residue 
$B$ is given by the {\it mapping cone} of the previous diagram
$$ 0\to G_1^\vee \to F_1^\vee\oplus G_2^\vee \to \dots\to F_{n-1}^\vee\oplus G_n^\vee$$
twisted by $-d_1-\dots-d_n$.

As a consequence, the $h$-vectors $Dh_A$, $Dh_B$, and $Dh_Z$ of sets $A,B$ linked by a complete intersection $Z$ as above are related by the following formula
\begin{equation}\label{hlink}
Dh_B(i) + Dh_A(d_1+\dots +d_n- n-i) = Dh_Z (i).
\end{equation}

The previous construction generalizes to the case in which $Z$ is arithmetically Gorenstein.

\begin{defn0} We say that a finite set $Z$ is \emph{arithmetically Gorenstein} if a minimal resolution of $I_Z$ is auto-dual, i.e. for all $i$ the dual of the map 
$F_i\to F_{i-1}$ is, up to twist, the map $F_{n-i+1}\to F_{n-i}$ (here we take $F_0=I_Z$). \\
In particular, we have that $F_n$ has rank $1$ and, for all $i$, $F_i$ is dual to $F_{n-i}$. This implies that the $h$-vector $Dh_Z$ 
of $Z$ is symmetric, i.e. 
$$Dh_Z(i)=Dh_Z(s-i)\quad\forall i.$$
where $s$ is the maximum such that $Dh_Z(s)>0$.
\end{defn0}

All complete intersection sets are arithmetically Gorenstein (but the converse is false).

If $Z=A\cup B$ is arithmetically Gorenstein, then the mapping cone of the diagram obtained by the resolutions of $I_Z,I_A$ provides, as in the case
in which $Z$ is a complete intersection, a resolution for $I_B$.

\begin{prop0} \label{DGO} (see \cite{DavisGerOre85})  If the $h$-vector $(h_0,h_1,\dots,h_s)$ of $Z$ is symmetric and $Z$ satisfies $CB(s)$,
then $Z$ is arithmetically Gorenstein.
\end{prop0} 

\subsection{Rank $2$ bundles and Serre construction}\label{sec:Serre}  When a set $Z$ of points lies in a smooth surface $S\subset\Pj^n$ the Serre
construction provides a link between the Cayley-Bacharach property and the existence of rank $2$ vector bundles on $S$ associated to $Z$.
We list in this section just the aspects of the connection that will be necessary for our analysis.

We recall that if $\Ec$ is a vector bundle of rank $2$ on a smooth surface $S$, with Chern classes $c_1(\Ec)\in Pic(S)$ and $c_2(\Ec)\in \mathbb Z$,
and $D$ is a divisor on $S$, then the Chern classes of the twist $\Ec(D)$ are given by
\begin{equation}\label{chern}
c_1(\Ec(D))=c_1(\Ec)+2D \qquad c_2(\Ec(D)) =c_2(\Ec)+ D\cdot c_1(\Ec) + D^2.
\end{equation}

\begin{prop0} \label{Serre} (Serre construction) Let $S$ be a smooth surface in $\Pj^n$, with canonical divisor $K$ and hyperplane divisor $H$. 
Let $Z\subset S$ be  a finite set which satisfies property $CB(d)$. Then there is a rank $2$ vector bundle $\Ec$ on $S$, with Chern classes
$c_1(\Ec)=dH+K$ and $c_2(\Ec)=\ell(Z)$ such that $Z$ is the zero-locus of a global section of $E$. The ideal sheaf $\Ic_{Z,S}$ of $Z$ in $S$ fits in an exact sequence
$$ 0 \to \Oc_S \to \Ec\to \Ic_{Z,S}(dH-K)\to 0,$$
in which the first map represents the global section that vanishes on $S$.\\
Conversely, if a finite set $Z$ is the zero-locus of a global section of a rank $2$ bundle $\Ec$ as above, with first Chern class $dD+K$, then $Z$ satisfies
the Cayley-Bacharach property $CB(d)$.
\end{prop0}

We will apply the previous proposition when $S$ in $\Pj^4$ is a complete intersection of two quadrics, so that $K=-H$. In this case the exact sequence reads:
\begin{equation}\label{seqSerre}
0 \to \Oc_S \to \Ec\to \Ic_{Z,S}((d+1)H)\to 0.
\end{equation}

For the proof, we refer to  \cite{Brun80}.

The existence of a section of a vector bundle $\Ec$ whose zero-locus is finite is regulated by the following proposition (see Remark 1.0.1 of \cite{Hart78}).

\begin{prop0} \label{codim2} Assume that $\Ec$ has non-trivial global sections. If all global sections
 of $\Ec$ vanish in  infinitely many points, then there exists an effective divisor $D$ on $S$ such that all sections of $\Ec$ are
 given by the product of a global section of $\Ec(D)$  times a global section of  $\Oc(D)$.
 \end{prop0}

\section{Forms of degree $4$ in five variables}\label{sec:quartics}

We turn now to the case $ n = 4, d = 4 $ and let $ T \in S^{4} \C^{5} $. Thus $T$ can be seen as a polynomial of degree $4$ in five variables, 
which is associated to a hypersurface of degree $4$ in $\Pj^4$.

\begin{rem0} For a general $ T \in S^{4} \C^{4} $, according to the Alexander-Hirschowitz Theorem \cite{AlexHir95}, the rank is $15$, while
the case of rank $14$ is defective, in the sense of \cite{CCi02a}. In both cases the dimension of the secant variety is smaller than the dimension of the corresponding
abstract secant variety, so that the identifiability cannot hold (this is a well known fact, see e.g. \cite{COtt12}, Proposition 2.2 and its proof).\\
On the other hand, if a specific $ T $ admits a non-redundant decomposition $ A $ of cardinality $ r \leq 8 $, then Theorem \ref{thm:kr} can 
be applied to establish the identifiability and rank of $ T $, while if $ r = 9 $, then one can refer to the criterion developed in \cite{AngeCVan18} 
and recalled in Theorem \ref{thm:krquartics}. 
\end{rem0}

Therefore we assume, from now on, that $ 10 \leq r \leq 13 $.

Let $ A = \{P_{1},\ldots,P_{r}\} \subset \Pj^{4} $ be a finite set that computes $ T $. We often assume that  $A$ satisfies the following conditions:
\begin{equation*}\label{**}
\qquad
\begin{cases} \begin{array}{ll}
(i) &  A \mbox{  is non-redundant; }\\
(ii) &  k_{1}(A) = 5 ; \\
(iii) &  k_2(A) = r.
\end{array}\end{cases}
\end{equation*}

It is a standard fact that when $ r\leq 21$, for $A$ in a Zariski open subset of $(\Pj^4)^r$, then $A$ satisfies the previous conditions.

Conditions $(i), (ii), (iii)$ imply that  the Hilbert function of $ A $ and its first difference  verify
\begin{equation}\label{eq:hDhA}
\begin{tabular}{c|ccccc}
$j$ & $0$ & $1$ & $2$ &   $3$ &   $\dots$ \\  \hline
$h_{A}(j)$ & $1$ & $5$ &   $r$ &   $r$ &   $\dots$ \cr
$Dh_{A}(j)$ & $1$ & $4$ &   $r-5$ &   $0$ & $\dots$ \cr
\end{tabular}.
\end{equation}

We notice that conditions $(i), (ii), (iii)$ can be easily controlled by a linear algebra algorithm (see the code in {\tt{ancillary.txt}} that can be implemented in Macaulay2, 
\cite{Macaulay2}). \\
When $A$ is sufficiently general, then conditions $(i), (ii), (iii)$ hold. Namely, a general set $A'$ of $r'$ points in $\Pj^4$, $5\leq r'\leq 15$
satisfies $h_{A'}(1)=5$, $h_{A'}(2)=r$. It is clear that if $A$ is general, then any subset of $A$ is general.
\smallskip

The following result is the basis for our analysis. It is a consequence of Theorem 1.2 of \cite{Ball19}. We provide here a different proof.

\begin{prop0} \label{quads} Let $A,B$ be non-redundant sets that compute $T$. Assume that $A$ satisfies conditions $(i), (ii), (iii)$, and assume that $A,B$ are disjoint. Then $B$ has length equal to $A$. \\
Let $\mathcal Q$ be the base locus of the linear system of quadric hypersurfaces passing through $A$. Then $\mathcal Q$ contains $B$ too. 
\end{prop0}
\begin{proof}
Set $\ell(A)=r$ and $ Z = A \cup B \subset \Pj^{4} $. Notice that, being $ A \subset Z $, from \eqref{eq:hDhA} we get that $ Dh_{Z}(1) = 4 $ and $ Dh_{Z}(2) \geq r-5 $. Moreover, 
Proposition \ref{d+1}  implies  $ Dh_{Z}(5) \geq 1 $. \\
Since $ A $ and $ B $ are non-redundant for $ T $, then by Proposition \ref{CBconseq}, $ Z $ has the Cayley-Bacharach property $ CB(4) $. 
Therefore from Theorem \ref{GKRext} we have  
$$ Dh_{Z}(3) + Dh_{Z}(4) + Dh_{Z}(5) \geq Dh_{Z}(0) + Dh_{Z}(1) + Dh_{Z}(2) \geq r. $$
Moreover, from the chain of inequalities
$$2r \geq \ell(Z) = \sum_{j \in \Z} Dh_{Z}(j) \geq \sum_{j=0}^{5} Dh_{Z}(j) \geq r + Dh_{Z}(3) + Dh_{Z}(4) + Dh_{Z}(5),$$
we deduce that
$$ Dh_{Z}(3) + Dh_{Z}(4) + Dh_{Z}(5) \leq r. $$
It turns out that 
$$ Dh_{Z}(3) + Dh_{Z}(4) + Dh_{Z}(5) = r. $$
Necessarily it has to be that $ Dh_{Z}(2) = Dh_{A}(2) = r-5 $, $ \ell(Z) = 2r $ and $ \ell(B) = r$. In particular, $ h_{Z}(2) = r$ and $ (I_{Z})_{2} = (I_{A})_{2} $. 
\end{proof}

Thus, if we can compute that $\mathcal Q$ is finite, of length $<2r$, we obtain the non-existence of $B$. This can be easily achieved by a computer-aided procedure,
as described in Remark \ref{generali}.

\subsection{Non-empty intersection}\label{nonO}

First we consider the case $A\cap B\neq \emptyset$, following the argument used several times in \cite{AngeC20}.
\smallskip

Assume that $ A = \{P_{1}, \ldots, P_{r}\}$ and let $B$ be another non-redundant decomposition of $T$ with $ s = \ell(B)\leq r$ and define $Z=A\cup B$. \\
Assume that  the intersection $A\cap B$ is not empty.\\
Then we can reorder the points of $A$ so that $B=\{P_1,\dots,P_j,P'_{j+1},\dots,P'_s\}$ with $j\geq 1$ 
and $P'_i\notin A$ for $i=j+1,\dots,s$. For any choice of representatives (i.e. coordinates) $T_1,\dots T_r$, and $ T'_{j+1},\dots,T'_{s}$ for the projective points 
$v_4(P_1),\dots,v_4(P_r)$, and $ v_4(P'_{j+1}),\dots, v_4(P'_s)$ respectively, there are non-zero scalars $a_i$'s, $b_i$'s such that
$$\begin{matrix}
T & = &a_1T_1+\dots+a_rT_r \\ T & = & b_1T_1+\dots+b_jT_j+b_{j+1}T'_{j+1}+\dots+b_sT'_s.
\end{matrix}$$
Define:
\begin{multline*} T_0=(a_1-b_1)T_1+\dots+(a_j-b_j)T_j+a_{j+1}T_{j+1}+\dots+a_rT_r
\\  = b_{j+1}T'_{j+1}+\dots+b_sT'_s.\end{multline*}
Now $T_0$ has the two decompositions $A$ and $B'=\{P'_{j+1},\dots, P'_s\}$, which are disjoint.  If $B'$ is redundant, then after 
rearranging the points, we may assume
$T_0=c_{j+1}T'_{j+1}+\dots+c_tT'_t$ for some  $t<s$, so that:
\begin{multline*} T=b_1T_1+\dots +b_jT_j+T_0=  b_1T_1+\dots +b_jT_j+c_{j+1}T'_{j+1}+\dots+c_tT'_t,
\end{multline*}
against the  fact that $B$ is non-redundant. Thus $B'$ must be non-redundant. \\
If $ A $ is non-redundant, define $ A' = A $. \\
If $A$ is redundant, since the points $v_4(P_1),\dots ,v_4(P_r)$ are linearly independent, then some coefficient $(a_i-b_i)$ is $0$. 
In this case, we may assume $(a_i-b_i)=0$ if and only if $i=1,\dots,q\leq j$, so we get a non-redundant decomposition $A'=\{P_{q+1},\dots,P_r\} $ of $ T_{0} $. \\
In conclusion, we find that $ T_{0} $ has two different non-redundant  decompositions $A', B'$, with $ A' \subset A $ and $ \ell(B') \leq \ell(A') $, $ \ell(B') < r $.
Notice that since $A'\subseteq A$, then conditions $(i), (ii), (iii)$ hold for $A'$. 

If $\ell(A')\leq 8$, then $T_0$ cannot exist, by Theorem \ref{thm:kr}. Thus $\ell(A')\geq 9$.

If  $\ell(B')<\ell(A')$ then  $T_0$ cannot exist, by Proposition \ref{quads}.
\smallskip

In conclusion we have the following

\begin{prop0} \label{inters}
If $T$ has two non disjoint decompositions $A,B$ of length $13\geq \ell(A)\geq\ell(B)$ and $A$ satisfies condition $(i), (ii), (iii)$, then $\ell(B)=\ell(A)\geq 9$. \\
Moreover, there are disjoint subsets $A'\subset A$ and $B'\subset B$, of the same length $r'\geq 9$, which are 
non-redundant decompositions of a form $T_0$ such that $T=T_0+ \sum_{i=0}^{r-r'} a_iv_4(P_i)$ where $\{P_1,\dots, P_{r-r'}\}=A\cap B$.
\end{prop0}

It follows that in many arguments, after replacing $T$ with $T_0$, we will assume that $A,B$ are disjoint.\\
Mixing Proposition \ref{quads} and Proposition \ref{inters} we get a first result which determines the rank of $T$. 

\begin{thm0}\label{rank} Let $ T \in S^{4} \C^{5} $ be a form with a decomposition $A$ of length $r\leq 13$, which satisfies conditions $(i), (ii), (iii)$ above.
Then $T$ cannot have a decomposition of length smaller than $r$. In other words, $r$ is the (Waring) rank of $T$.
\end{thm0}

Next, we turn the attention to the identifiability of $T$.

\subsection{Case  r = 9, 10, 11} 

With the previous notation, we assume that there exists a second, non-redundant decomposition $B$ of $T$ of length $r$.

\begin{rem0}\label{generali} Let $\mathcal Q'$ be the base locus of the linear system of quadric hypersurfaces passing through a finite set $A'$.
Let $A'\subset \Pj^4$ be general with $\ell(A') = 9$. Then it is a consequence of the Minimal Resolution Conjecture (which holds in 
$\Pj^4$, see \cite{Walter95}) that the homogeneous ideal $I_{A'}$ is generated by quadrics.
Indeed for a general set $A$ of $9$ points, $I_A$ has a free resolution of the form
\begin{equation}\label{eq:idA9}
0 \rightarrow R(-6)^{\oplus 4} \longrightarrow R(-5)^{\oplus 12} \longrightarrow R(-4)^{\oplus 9} \oplus R(-3)^{\oplus 4} \longrightarrow R(-2)^{\oplus 6} \longrightarrow I_{A'} \rightarrow 0.
\end{equation}
In this case,  $\mathcal Q'$ coincides with $A'$ (see the file {\tt{ancillary.txt}})

Similarly, for a general set $A'$ of $10$ points,  the homogeneous ideal $I_{A'}$ is generated by quadrics and it has a free resolution of the form
\begin{equation}\label{eq:idA10}
0 \rightarrow R(-6)^{\oplus 5} \longrightarrow R(-5)^{\oplus 16} \longrightarrow R(-4)^{\oplus 15} \longrightarrow R(-2)^{\oplus 5} \longrightarrow I_{A'} \rightarrow 0.
\end{equation}
Thus also in this case $\mathcal Q'$ coincides with $A'$ (see the file {\tt{ancillary.txt}})

For a general set $A'$ of $11$ points, the ideal $I_{A'}$ has a free resolution of the form
\begin{equation}\label{eq:idA11}
0 \rightarrow R(-6)^{\oplus 6} \longrightarrow R(-5)^{\oplus 20} \longrightarrow R(-4)^{\oplus 21} \longrightarrow R(-3)^{\oplus 4} \oplus R(-2)^{\oplus 4} 
\longrightarrow I_{A'} \rightarrow 0.
\end{equation}
In this case $\mathcal Q'$ does not coincide with $A'$. On the other hand, for a general choice of $A'$ one computes (see the file {\tt{ancillary.txt}}) that  four general quadrics
in $I_{A'}$ intersect properly, so that $\mathcal Q'$ is a finite set of length $16$. 
\end{rem0}

As a consequence, we can use the second part of proposition \ref{quads} to prove the identifiability of $T$. We add the following condition for $A$.
\smallskip

$(iv)$ \hskip.5cm  the base locus $\mathcal Q$ of the system of quadrics containing $A$ is finite.
\smallskip

Notice that a general set of $r\leq 11$ points will satisfy condition $(iv)$. Moreover, if $A$ satisfies condition $(iv)$, then trivially all subsets of $A$
also satisfy the condition.

\begin{prop0}\label{allscalars} Let $T$ be a form of degree $4$ in $5$ variables. Assume:
\begin{equation} \label{eq:idC}  T  = a_1T_1+\dots+a_rT_r \end{equation}
where $T_i$'s are powers of linear forms, corresponding to points $v_4(P_1),\dots,v_4(P_r)$.\\
 Assume $r\leq 11$ and
assume that the set $A=\{P_1,\dots,P_r\}$ satisfies conditions $(i), (ii), (iii), (iv)$.\\
Then $T$ is identifiable, of rank $r$. I.e., \eqref{eq:idC} provides the unique decomposition of length $r$ of $T$ (up to multiplication by a scalar and rearranging).
\end{prop0}
\begin{proof} The assumptions say that $A$ is a non-redundant decomposition of $T$.\\ 
By condition $(iv)$, four general quadrics containing $A$ meet in a finite set, so the base locus $\mathcal Q$ of the system of quadrics containing $A$ is a finite set
 of length at most $16$.  \\
Consider a  second decomposition $B$ of $T$ of length $\leq r$. Since $T$ has rank $r$ by Theorem \ref{rank}, then in fact $\ell(B)=r$.\\
 If $A,B$ are disjoint, then $A\cup B$ is contained in the intersection $\mathcal Q$ of the quadrics containing $A$, by Proposition \ref{quads}, 
 so that $\ell(B)\leq 16-\ell(A)$ and  we get a contradiction. \\
 If $A,B$ are not disjoints, by Proposition \ref{inters} there are disjoint subsets $A'\subset A$ and $B'\subset B$, of the same length $r'\geq 9$,
 which are both decompositions of a form $T_0$. Since $A'$ also satisfies conditions $(i), (ii), (iii), (iv)$, then from the previous argument we get
 that  $B'$ is contained in the intersection $\mathcal Q'$ of quadrics containing $A'$, which is a set
 of length at most $16-\ell(A')\leq 7$, a numerical contradiction.
 \end{proof}

 Notice that the identifiability of $T$ follows once we test that $A$ is sufficiently general to satisfy conditions $(i), (ii), (iii), (iv)$.  
  The file  {\tt ancillary.txt} contains an algorithm that checks if a specific $A$ of length $r\leq 11$ satisfies the four conditions. 
  \smallskip
  
  It is clear that $(ii), (iii), (iv)$ only deal with the geometry of $A$. It follows that if $A$ satisfies $(ii), (iii), (iv)$, then the identifiability
  holds for all tensors in the span of $v_4(A)$, for which $A$ is non-redundant, i.e. for all $T$ of the form
 $ T  = a_1T_1+\dots+a_rT_r $, where all the coefficients $a_i$'s are non-zero.

\subsection{Case  r = 12}

When $r=12$, we still know from Theorem \ref{rank} that when $A$ satisfies $(i), (ii), (iii)$ then $T$ has rank $12$.

The situation for the identifiability of $T$ changes. Indeed the ideal $I_A$ of a general set of $12$ points in $\Pj^4$ has free resolution of the form 
\begin{equation}\label{eq:idD12}
0 \rightarrow R(-6)^{\oplus 7} \longrightarrow R(-5)^{\oplus 24} \longrightarrow R(-4)^{\oplus 27} \longrightarrow R(-3)^{\oplus 8} \oplus R(-2)^{\oplus 3} 
\longrightarrow I_{A} \rightarrow 0.
\end{equation}
and the Hilbert function and its difference are given by:
\begin{equation}\label{eq:hDhE}
\begin{tabular}{c|ccccc}
$j$ & $0$ & $1$ & $2$ &   $3$ &   $\dots$ \\  \hline
$h_{A}(j)$ & $1$ & $5$ &   $12$ &   $12$ &   $\dots$ \cr
$Dh_{A}(j)$ & $1$ & $4$ &   $7$ &   $0$ & $\dots$ \cr
\end{tabular}.
\end{equation}

Thus there are only three independent quadrics in $I_A$, therefore the intersection $\mathcal Q$ contains (and in general it is equal to) a curve of degree $8$.
It follows that  condition $(iv)$ cannot hold, and we cannot use the argument above to exclude the existence of a second decomposition $B$ of length $12$.\smallskip

We replace condition $(iv)$ with
\smallskip

$(iv')$ \hskip.5cm  for all subsets $A'\subset A$ of length at most $11$, the base locus $\mathcal Q'$ of the system of quadrics containing $A'$ is finite.
\smallskip

Notice that $(iv')$ can be checked only considering subsets of length $11$.\smallskip

When we start with  a decomposition $A$ of length  $r=12$, testing the identifiability the situation becomes much more complicated.

By \cite{COttVan17a}, we know that a decomposition of a general form of degree $4$ in $5$ variables and rank $12$ is unique. But in order to get
a criterion which determines the uniqueness of the decomposition, it is not sufficient to look at the points $P_i$'s. We will see that for a general
choice of the set $A=\{P_1,\dots,P_{12}\}$  there are points $T$ in the span of $v_4(P_1),\dots, v_4(P_{12})$ for which $A$ is the unique
decomposition, and points $T'$ for which $A$ is non-redundant, but not unique.

\begin{exa0} \label{fond12} Let $A=\{P_1,\dots,P_{12}\}$ be a general set of $12$ points in $\Pj^4$.\\
The exact sequence \eqref{eq:idD12} provides a minimal resolution for the ideal of $A$, which is thus generated in degree $3$.
This means that if we take $3$ independent quadrics and a general cubic $F$ in $I_A$, the intersection $Z$ is a finite set of $24$ distinct points in $\Pj^4$.
$Z$ is a complete intersection set of points. The residue $B$ is a set of $12$ points, different from $A$. Indeed, since $Z$ is reduced, 
for a general choice of $A$ and the cubic $F$, the two sets $A,B$ are disjoint.\\
We claim that the spans of $ v_4(A)$ and $v_4(B)$ meet in a point $T$ respresenting a form of degree $4$ in $5$ variables, with rank $12$
for which $A,B$ are two different minimal (non-redundant) decompositions. In other words, $T$ is non-identifiable.

In order to prove the claim, thanks to Theorem \ref{rank} and Proposition 2.19 of \cite{AngeC20}, 
it is sufficient to prove $v_4(A\cup B)=v_4(Z)$ is not linearly independent, i.e. $h_Z(4)<24$.
This is well known: it follows immediately from the free resolution of a complete intersection of type $2,2,2,3$ (i.e. a complete intersection of $3$ quadrics and one cubic) in $\Pj^4$,
which is given by the Koszul complex: \small{
$$
0 \rightarrow R(-9) \rightarrow R(-7)^{\oplus 3}\oplus R(-6) \rightarrow R(-5)^{\oplus 3}\oplus R(-4)^{\oplus 3}\rightarrow R(-3) \oplus R(-2)^{\oplus 3} \rightarrow I_{Z} \rightarrow 0.
$$
}
\end{exa0}

We want to prove that, at least when $A$ is general, the previous example provides the unique possibility for a form
$T\in\langle v_4(A)\rangle$, for which $A$ is non-redundant, to be non-identifiable.
\smallskip

\begin{rem0}\label{disj} 
Assume that $A$ satisfies conditions $(i), (ii), (iii), (iv')$.
If $T$ has a second decomposition $B$ of length $12$, then by arguing as in Proposition \ref{allscalars}
we can prove that $A\cap B$ is empty.

Indeed if $A\cap B$ is non-empty, then  by Proposition \ref{inters} there are disjoint subsets $A'\subset A$ and $B'\subset B$, of the same length $r'\leq 11$,
 which are both decompositions of a form $T_0$. Since $A$ satisfies $(i), (ii), (iii), (iv')$, then $A'$ also satisfies conditions $(i), (ii), (iii), (iv)$, 
 so that $B'$ is contained in the intersection the intersection $\mathcal Q'$ of quadrics containing $A'$, which is a set
 of length at most $16-\ell(A')\leq 7$, a contradiction.
\end{rem0}

\begin{prop0} If $A$ is general, then the intersection of the quadric hypersurfaces  in $I_A$ is an irreducible
(complete intersection) curve $C$ of degree $8$ and arithmetic genus $5$.
\end{prop0}
\begin{proof} The irreducibility of $C$ follows by a standard application of Bertini's Theorem, see \cite{Hartshorne} Lemma V.1.2.
\end{proof}

For the rest of the section, we will always assume that the decomposition $A$ satisfies conditions $(i), (ii), (iii), (iv')$ and moreover
\smallskip

 $(v)$ \hskip.1cm   the base locus  of the system of quadrics through $A$ is an irreducible curve, 

    \hskip 0.7cm and the  homogeneous ideal $I_A$ has a resolution as \eqref{eq:idD12}.
\smallskip

Condition $(v)$ can be checked with the aid of computer algebra packages (see the file {\tt ancillary.txt} for details).

\begin{prop0}\label{linkedsch} Assume that $A$ satisfies conditions $(i), (ii), (iii), (iv'), (v)$. If $T$ has a second decomposition $B$ of length $12$,
then $A\cup B$ is a complete intersection of type $2,2,2,3$.
\end{prop0}
\begin{proof} We know that $A\cup B=Z$ is a set of length $24$ which satisfies the Cayley-Bacharach property $CB(4)$. Thus the 
difference of the Hilbert function of $Z$ satisfies:
$$\begin{array}{ll}
Dh_Z(0)=1,  \quad Dh_Z(1)=4, \quad  Dh_Z(2)\geq 7 \\ 
Dh_Z(3) + Dh_Z(4)+ Dh_Z(5)  \geq  Dh_Z(0)+Dh_Z(1)+ Dh_Z(2)  \geq 12, \\ 
Dh_Z(4)+ Dh_Z(5) \geq  Dh_Z(0)+Dh_Z(1)  = 5,\\
\sum_{i=0}^5 Dh_Z(i)  = 24.
\end{array}$$
which imply that $Dh_Z(2)=Dh_A(2)=7$ and also  $h_Z(3)\leq19$. Thus the ideal $I_Z$ coincides with $I_A$ up to degree $2$, hence
$I_Z$ contains three independent quadrics, whose intersection is an irreducible curve $C$. Moreover $I_Z$ contains a cubic $F$
which is not spanned by the quadrics of $I_Z$. Hence $C\cap F$ is a set of $24$ points, so that $Z=F\cap C$.
\end{proof}

The linkage determines a resolution of the ideal $I_B$ of $B$, via the mapping cone procedure  (see Section \ref{sec:liaison}).
 Indeed, from the resolution of $A$ we get a diagram:  
{\small $$
\begin{matrix}
0 &\to& R(-9) &\to& \begin{matrix}R(-7)^{\oplus 3} \\ \oplus \\ R(-6) \end{matrix}&\to& 
\begin{matrix} R(-5)^{\oplus 3}\\ \oplus \\ R(-4)^{\oplus 3} \end{matrix} &\to& 
\begin{matrix} R(-3) \\ \oplus \\ R(-2)^{\oplus 3}\end{matrix}  &\to& I_Z &\to& 0 \\ 

 & & \big\downarrow & & \big\downarrow & &  \big\downarrow & & \big\downarrow & & \big\downarrow  & &  \\

0 &\to& R(-6)^{\oplus 7} &\to& R(-5)^{\oplus 24} &\to& 
R(-4)^{\oplus 27} &\to& \begin{matrix}R(-3)^{\oplus 8}\\ \oplus \\ R(-2)^{\oplus 3}\end{matrix}  &\to& I_{A} &\to& 0
\end{matrix}$$}
\noindent where $Z=C\cap F$. We get that a resolution of $I_B$ is numerically identical to the resolution of $I_A$.

\begin{exa0} \label{nonredund12} Fix a general set $A$ of $12$ points in $\Pj^4$. $A$ is contained in three independent quadrics, 
whose intersection is a general (smooth) canonical curve $C$. Cubics containing $A$ and not $C$ determine a $7$-dimensional projective space 
$\Pj((I_A)_3/(I_C)_3)$. From the coordinates of the set $A$ one computes the ideal $I_A$. From the choice of a cubic hyperpsurface $F$
through $A$ and the mapping cone procedure, one computes the ideal $I_B$ of the set of $12$ points, residue of the intersection $F\cap C$ with respect to $A$.
If everything is general enough, then $B\cap A$ is empty and $h_Z(4)=23$, where $Z=A\cup B$. 
It follows that in the polynomial ring $R=\C [x_0,\dots,x_4]$ the subspaces $(I_A)_4$ and $(I_B)_4$ of $R_4$ are both $58$-dimensional,
but their sum $U=(I_A)_4+(I_B)_4$ is a linear subspace of codimension $1$ in $R_4$. Coefficients for the equation of $U$
provide coefficients for the form $T$ in the intersection of $\langle v_4(A)\rangle$ and $\langle v_4(B)\rangle$. We refer to the file {\tt ancillary.txt} for details.\\
Thus, from the coordinates of $A$ and $B$ one determines $T$ explicitly, and also an expression of $T$ as combination
of rank $1$ tensors representing the points of $A$. Experiments prove the existence of $T\in\langle v_4(A)\rangle\cap \langle v_4(B)\rangle$ for which no coefficients
in the decomposition given by $A$ vanish. We refer for that to the file {\tt ancillary.txt}.\\
If $A$ is general enough, then the span of tangent spaces to the Veronese variety $v_4(\Pj^4)$ at the points of $v_4(A)$ determines 
a linear subspace of the expected projective dimension $59$  in the projective space spanned by $v_4(\Pj^4)$. This is easy to check
directly for the points $A$ above (see the  file {\tt ancillary.txt}). From Remark \ref{stacked}
 it follows  that there are no positive-dimensional  families  of decompositions for the form $T$.
\end{exa0} 

\begin{coro0} \label{cor10} Fix a set $A$ of $12$ points in $\Pj^4$, which satisfies condition $(i), (ii),$ $ (iii), (iv'), (v)$. 
Then the span of $v_4(A)$, which is a $\Pj^{11}$, contains an
irreducible  $7$-dimensional family $\Tc$ of $T$ for  which $A$ is a non-redundant decomposition, and such that
they  have a second decomposition $B$ obtained as in Example \ref{nonredund12}. 
\end{coro0}
\begin{proof} The family $\Bc$ of sets $B$ of $12$ points linked to $A$ in a complete intersection of type $2,2,2,3$ is irreducible and $7$-dimensional.
For a general $B$ in the family, $h_{A\cup B}(4)=23$, so that $\langle v_4(A)\rangle$ and $\langle v_4(B)\rangle$ meet in one point $T$.
Call $\Tc$ the family of points $T$ that we obtain in this way.\\
In specific examples one computes that  the tangent spaces to $v_4(\Pj^4)$ at the points of $v_4(A)$ span a space of dimension $59$
(see the  file {\tt ancillary.txt}). Thus, there are no positive dimensional families of decompositions $B$ which map to the same general $T$
in the family.
\end{proof}

We describe an algorithm that tests if the decomposition $A$ satisfies conditions $(i), (ii), (iii), (iv'), (v)$ 
and moreover if the given form $T$ belongs to the family $\Tc$. The algorithm has been implemented in the file {\tt ancillary.txt}.\\
Since all the conditions controlled by the algorithm are Zariski open, then the algorithm is {\it effective}, in the sense of \cite{COttVan17b}:
it will give a positive answer for all $T$ lying in a Zariski open, dense  subset of the space of all quartics in five variables.
\medskip

{\bf Algorithm} 
\begin{itemize}
\item[input:] a Waring expression $T = \lambda_{1}L_1^4+\dots+\lambda_{12}L_{12}^4$ of a quartic $T$ in five variables, 
where each $L_i$ is a linear form, represented by the vector $u_i$ of coefficients, and the corresponding set $A\subset \Pj^2$. 
\item[1]  {\it test (i)}: check that $\dim \langle v_4(u_1),\dots, v_4(u_{12})\rangle = 12$;
\item[2] {\it test (ii)}: check that $k_1(A) = 5$;
\item[3] {\it test (iii)}: check that $k_2(A) = 12$;
\item[4] {\it Terracini test}: check that the span of the tangent spaces of $v_4(\Pj^2)$ at the points of $v_4(A)$ has the expected dimension $59$;
\item[5] {\it test (iv')}: check that, for all $i$, the base locus of the system of quadrics through the points corresponding to  
$L_j$, $j\neq i$, is finite, of length $\leq 16$;
\item[6] {\it test (v)}: check that the base locus of the system of quadrics through $v_4(A)$ is an irreducible curve;
\item[7] {\it parametrization of $\Tc$}: construct a parametrization of alternative decompositions $B$ of forms in the span of $v_4(A)$;
\item[8] {\it final test}: test that for no choice of the parameters for $B$ the form $T$ is orthogonal to generators of $(I_B+I_A)_4$.
\end{itemize}
If any of the tests above answers negatively, then the algorithm terminates by returning that it cannot handle the expression  $T = \lambda_{1}L_1^4+\dots+\lambda_{12}L_{12}^4$.\\
Otherwise the algorithm returns that the expression is minimal and unique.

{\bf end of Algorithm}
\medskip

Notice that steps $1-4$ and $8$ require only standard methods of linear algebra for matrix analysis, while steps $5,6$ need 
also computer algebra packages.  \\
Notice also that, when the tests $1-6$ are positive, then a negative answer to test $8$ can be refined to provide generators for the
ideal of a second decomposition $B$ of length $12$ for $T$.

\begin{exa0} {\it An identifiable quartic}\\
Let 
$$ T_{1} = \sum_{i=1}^{12} L^{4}_{i} =  $$
{\small{$$ = [2454x^4-14837x^3y-9546x^2y^2+12272xy^3+5779y^4+9852x^3z+11840x^2yz+4479xy^2z-6699y^3z+$$
$$ + 5245x^2z^2+7347xyz^2+979y^2z^2-14274xz^3-10753yz^3-15547z^4+3625x^3t+1511x^2yt-7021xy^2t+$$
$$+8756y^3t-12116x^2zt-11133xyzt-4526y^2zt-8491xz^2t+12057yz^2t-9401z^3t-10613x^2t^2+$$
$$-6878xyt^2+8208y^2t^2+3405xzt^2+10766yzt^2-13732z^2t^2+14028xt^3-9572yt^3-11158zt^3-2774t^4+$$
$$-5103x^3w+5136x^2yw+10632xy^2w-15393y^3w-4914x^2zw+8047xyzw-4020y^2zw-1609xz^2w+$$
$$+14390yz^2w-5791z^3w+8743x^2tw-14600xytw+11388y^2tw+6681xztw+15846yztw+9266z^2tw+$$
$$+3649xt^2w-4887yt^2w+12361zt^2w+8699t^3w+12211x^2w^2-10563xyw^2-13952y^2w^2+$$
$$+2139xzw^2-12182yzw^2-7237z^2w^2-113xtw^2-1224ytw^2-2612ztw^2+13999t^2w^2-6977xw^3+$$
$$-8368yw^3+1738zw^3-14977tw^3+3637w^4]\quad\quad\quad\quad\quad\quad\quad\quad\quad\quad\quad\quad\quad\quad\quad\quad\quad\quad\quad\quad$$}}

\noindent be a quartic in five variables, where $ L_{1}, \ldots, L_{14} $ are linear forms represented by a random collection of 12 vectors generated in Macaulay2: 
$$ A = \left[\bf{u}_i \right]_{i=1}^{12}= {
\begin{bmatrix}
-1960 & 7185 &  2948 & 1986 &  -7270\\
8416 &  -14232 & 8567 &  14988 & -12297 \\
 4210  & -11055 & -6249 & 530  &  6066  \\
-6981 & 1313  & 6692 & 12883 & 4597  \\
8211 &  -5857  & 6853 & -5758 & -1890\\
8633 &  6895  & 14963 & 14147 & -405  \\
12697 & -10281 & 10647 & 1414  & 11296\\
-15107 & 4696 &  -6212 & 6064  & 8777\\
-14194 & -13431 & -2768 & 6063 &  -1066 \\
-687  & 7327  & 9904 &  11696 & 10323\\
-262  & -14530 & 5673 & 10210 & 5157\\
 -5397 & 6232  & -7867 & -10827 & -653
\end{bmatrix}.} $$
The matrix obtained by applying to each of the $ \bf{u}_i $'s the Veronese map $ v_{4} $ of degree $ 4 $ has full rank $ 12 $, so that test (i) has positive answer. Moreover, all the $ 792 $ maximal minors of $ A $ do not vanish, which implies that $ k_{1}(A) = 5 $. Similarly, by applying the Veronese map $ v_{2} $ to the rows of $ A $, we get a matrix with all the $ 455 $ maximal minors different from $ 0 $, i.e. test (iii) is successful. Terracini test provides a positive answer since the $ 60 \times 70 $ matrix associated to the span of the tangent spaces of $ v_{4}(\Pj^2) $ at the points of $ v_{4}(A) $ is of maximal rank. By removing one vector at a time from the rows of $ A $, we get that for each one the $ 12 $ sets of $ 11 $ points the base locus of the system of quadrics through them is finite dimensional, that is test (iv') is successful. Since the ideal of the base locus of the system of quadrics through the points represented by the $ \bf{u}_i $'s is irreducible, also test (v) answers positively. Final test translates into check that a linear system with $11$ equations and $ 8 $ unknowns (which represent the choices of affine parameters for alternative decomposition $ B $ for $ T_{1} $) admits only the trivial solution: since for $ T_{1} $ the matrix (MatEqns, in the notation of file {\tt ancillary.txt}) associated to this system is of full rank $ 8 $, then final test provides a positive answer for $ T_{1} $.  Therefore $ T_{1} $ is a rank-$ 12 $ identifiable quartic.
\end{exa0}

\begin{exa0}\label{nonid} {\it A non-identifiable quartic}\\
Consider the $12$ linear forms $L_1,\dots,L_{12}$ in $ x,y,z,t,w $, whose coefficients appear in the rows of the following matrix, generated randomly in Macaulay2:
$$ A = \left[\bf{u}_i \right]_{i=1}^{12}= {
\begin{bmatrix}
39  & -33 & -11 & 5 &  24\\
-30 & -28 & 44 & -19 & -32 \\
 -15 & 19 & -50 & 43 & 48 \\
-35 & 5  & -9 & -31 & 4 \\
30 & 8 &  5 &  -28 & -44 \\
-31& -33 &39 & 41 & 47  \\
-6 & -20 & 34 & -7 & 27 \\
-23 & -9 & 31 & -24 & 13\\
-16 & -34 & -36 & 4 &  5  \\
9  & 49 & 42  &-36 & 34 \\
14 & 50 & 44 & -4  &-17 \\
12 & 48 & 6  & -32 &-26 
\end{bmatrix}.} $$ 
\noindent Then, by reversing the procedure illustrated in the section, we constructed a quartic 
$$ T = \sum_{i=1}^{12} \lambda_{i}L^{4}_{i} $$
\noindent which has a second decomposition of length $12$. \\
The computer guided construction generated coefficients $ \lambda_{i} $'s which are considerably too long
to be reported here. We list them explicitly, together with the coefficients of $ T $ in the standard monomial basis of degree 4,
in the file {\tt{ancillary1.txt}}.\\
We checked that all the steps $1-6$ in the algorithm provide positive answers for $ T $. On the other hand, the final test does not succeed, being the $ 11 \times 8 $ matrix  MatEqns, constructed as  in the previous example, of rank $ 7 $. Indeed, $T$ has a second decomposition given by linear forms described in the following matrix: $ B = \left[\bf{v}_i \right]_{i=1}^{12}= $

{\tiny{
 $$\begin{bmatrix}
-7.35766+11.8909i & -6.51802+2.18767i &  11.0495+0.44639i &  -3.19976-2.02225i & 1 \\
-7.35766-11.8909i & -6.51802-2.18767i  & 11.0495-0.44639i  & -3.19976+2.02225i & 1 \\
 -2.99409+5.9125i & 3.81836+6.94562i  &  -1.56734+9.41357i & -7.62593-1.87791i &1 \\
-2.99409-5.9125i & 3.81836-6.94562i  &  -1.56734-9.41357i & -7.62593+1.87791i  & 1 \\
 .301069+0.983283i & 0.710226+.540009i  &  -1.10161-1.42291i & .489521+.792445i &  1 \\
.301069-0.983283i &  0.710226-.540009i  &  -1.10161+1.42291i & .489521-.792445i & 1 \\
 -1.65237  &  -.483634   &          1.8944     &        -1.66592     &      1 \\
 -0.430153   &        -1.88032      &       -0.131641    &       1.33533     &       1 \\
-.649097      &     -0.732771       &      0.791558     &       0.858328      &      1 \\
 -0.59154+.0315252i  & -0.891556+.00946299i  & 0.202355-0.0971274i  & 0.511408-0.174549i & 1 \\
 -0.59154-0.0315252i  & -0.891556-0.00946299i  & 0.202355+0.0971274i  & 0.511408+0.174549i  &1 \\
-0.512054       &    -2.11703      &       -1.36293    &       0.424299        &    1 
\end{bmatrix}. $$ 
}}

\noindent The coordinates of these points have been found by matching eigenvalues corresponding to the same eigenvector of certain companion matrices associated to the polynomial system given by a set of minimal generators of $ I_{B} $.
\end{exa0} 

We notice that the form $T$ constructed in Example \ref{nonid} is nevertheless identifiable over the reals.\\
Indeed, we know from Remark \ref{disj} and  Proposition \ref{linkedsch} that for any decomposition $B'$ of $T$ alternative to $A$, the union $A\cup B'$
is a complete intersection of three quadrics and one cubic. The quadrics are uniquely determined by $A$, while the cubic depends
on parameters. Our algorithm in the ancillary file {\tt ancillary.txt} for the final test
(step $8$ above) shows that there exists a unique choice of the parameters for the cubic that produce a set $B'$
such that $T$ is orthogonal to the to generators of $(I_{B'}+I_A)_4$.\\
Thus $A$ and the set $B$ described above are the unique  decompositions of length $12$ of $T$, and clearly $A$ is the unique one defined over the reals.\\
Compare with \cite{AngeBocciC18} for similar examples of tensors which are identifiable over the reals, but not identifiable over $\C$.

\subsection{Case  r = 13}

This turns out to be the most difficult case, because even when $A$ is general there are some degenerate cases to consider.
From now on, in this section we assume that $A$ satisfies conditions $(i), (ii), (iii), (iv')$ above, 
so that its Hilbert function is described in table \eqref{eq:hDhA}.

When $r=13$ the ideal $I_A$ of a general set of $r$ points in $\Pj^4$ has free resolution of the form
\begin{equation}\label{eq:idD}
0 \rightarrow R(-6)^{\oplus 8} \longrightarrow R(-5)^{\oplus 28} \longrightarrow R(-4)^{\oplus 33} \longrightarrow R(-3)^{\oplus 12} \oplus R(-2)^{\oplus 2} 
 \longrightarrow I_{A} \rightarrow 0.
\end{equation}

\begin{rem0}\label{esse}
For a general choice of $A$ there are only two independent quadrics containing $A$. Moreover the quadrics
 intersect in a smooth surface $S$ of degree $4$. \\
The surface $S$ has a well known structure. It corresponds to the blow up of $\Pj^2$ at $5$ points $Q_1,\dots,Q_5$, 
embedded with the linear system of cubics  through the points $Q_i$'s (a Del Pezzo surface).\\
The Picard group of $S$ is generated by the strict transform of a line $h$ and the five exceptional divisors $e_1,\dots,e_5$.
The hyperplane divisor is $H=3h-\sum e_i$. The canonical class of $S$ is $K_S=-H$.
\end{rem0}

\begin{exa0} \label{nodis} For $r=13$, even if the set $A$ is general, the span $L$ of $v_4(A)$ contains forms of rank $13$ with another
decomposition $B$ which intersects $A$.\\
Namely, write $A=\{P_1,\dots,P_{13}\}$ and consider $A'=\{P_1,\dots,P_{12}\}\subset A$. $A'$ is a general
set of $12$ points, and we saw in the previous section that the span of $v_4(A')$, which is a hyperplane in $L$,
contains forms $T'$ of rank $12$ with a second decomposition $B'$ of length $12$. If $T$ is a general form in the
line joining $T'$ with $v_4(P_{13})$, then necessarily $T$ has rank $13$ and two minimal decompositions
$A=A'\cup\{P_{13}\}$ and $B=B'\cup\{P_{13}\}$.\\
Thus, in order to prove that the decomposition $A$ is unique, it is necessary to compute, for $i=1,\dots,13$,
 the tensor $T'_i$ in the intersection of the line joining $v_4(P_i),T$ with the linear span of the $v_4(P_j)$'s, $j\neq i$, and
run the procedure illustrated in the previous section, to guarantee that $T'_i$ has no other decompositions of length $\leq 12$,
different from $A\setminus\{P_i\}$.
\end{exa0}

We want to prove that Example \ref{nodis} is the unique case in which $T$ has a second decomposition $B$, different from $A$
and not disjoint from $A$.

\begin{prop0} \label{rank>11} Let $A\subset\Pj^4$ be a set of $13$ points satisfying conditions $(i), (ii)$, $(iii), (iv')$ above. 
Then the span of $v_4(A)$ contains no tensors
of rank $\leq 11$, outside the $10$ dimensional subspaces generated by the images subsets of $A$ of length $11$.
\end{prop0}
\begin{proof} Assume that a form $T$, which sits in the span of $v_4(A)$ and not in the span of subsets of length
$11$ in $v_4(A)$,  has a second decomposition $B$ of length $\leq 11$. By proposition \ref{inters} we get $A\cap B=\emptyset$.
Then $Z=A\cup B$ satisfies the Cayley-Bacharach property $CB(5)$. Thus we have:
$$\begin{array}{ll}
\sum_{i=0}^5 Dh_Z(i)  \leq 24,  \\
Dh_Z(0)=1,  \quad Dh_Z(1)=4, \quad  Dh_Z(2)\geq 8 \\ 
Dh_Z(3) + Dh_Z(4)+ Dh_Z(5)  \geq  Dh_Z(0)+Dh_Z(1)+ Dh_Z(2)  \geq 13, \\ 
\end{array} $$
which all together provide a contradiction.
\end{proof}

\begin{prop0} \label{capnon0} Let $A\subset\Pj^4$ be a set of $13$ points satisfying conditions $(i), (ii)$, $(iii), (iv')$
above, and let $T$ be a form in the span
of $v_4(A)$, such that $A$ is non-redundant for $T$. Assume that 
$T$ has a second decomposition $B$ of length $s\leq 13$ such that $A\cap B\neq \emptyset.$\\
Then $s=13$ and there exists a point $P\in A$ and a form $T'$ in the span of $v_4(A\setminus\{P\})$ such that $T'$ has rank $12$
and $2$ disjoint decompositions of length $12$: $A\setminus\{P\}$ and $B'$, moreover $B=B'\cup\{P\}$.\\
In particular, $\ell(B)=13$ and $A\cap B$ is a singleton.
\end{prop0}
\begin{proof} The fact that $s=13$ follows from Proposition \ref{inters}.
With the usual notation, write $B=\{P_1,\dots,P_j,P'_{j+1},\dots,P'_{13}\}$, $j\geq 1$ and choose  
representatives (i.e. coordinates) $T_1,\dots T_{13}, T'_{j+1},\dots,T'_{s}$ for the projective points 
$v_4(P_1),\dots,v_4(P_{13})$ and $v_4(P'_{j+1}),\dots,$ $ v_4(P'_s)$ respectively. Then:
$$\begin{matrix}
T & = &a_1T_1+\dots+a_{13}T_{13} \\ T & = & b_1T_1+\dots+b_jT_j+b_{j+1}T'_{j+1}+\dots+b_sT'_s.
\end{matrix}$$
The form
\begin{multline*} T_0=(a_1-b_1)T_1+\dots+(a_j-b_j)T_j+a_{j+1}T_{j+1}+\dots+a_{13}T_{13}
\\  = b_{j+1}T'_{j+1}+\dots+b_sT'_s.\end{multline*}
 has the two decompositions $A$ and $B'=\{P'_{j+1},\dots, P'_s\}$, which are disjoint.  \\
If $\ell(A')\leq 11$, we get a contradiction
 with Proposition \ref{allscalars}. Thus 
 $\ell(A')=\ell(B')=12$, so that $j=1$ and the claim is proved.
\end{proof}

Now, let us turn to the case of disjoint decompositions $A$, $B$, with $A$ general of length $13$ and $B$ of length $\leq 13$.

\begin{prop0} \label{hilbdisj} Let $A\subset\Pj^4$ be a set of $13$ points satisfying conditions $(i), (ii)$, $(iii), (iv')$, and let $T$ be a form in the span
of $v_4(A)$, such that $A$ is non-redundant for $T$, and assume that $T$ has a second decomposition $B$ of length $s\leq 13$ such that 
$A\cap B= \emptyset.$ \\
Then $\ell(B)=13$, the ideal
of $Z=A\cup B$ coincides with the ideal of $A$ up to degree $2$, and $h_Z(3)\leq 21$.
\end{prop0}
\begin{proof} 
The set $Z$ satisfies the Cayley-Bacharach property $CB(5)$. Thus we have:
$$\begin{array}{ll}
\sum_{i=0}^5 Dh_Z(i)  \leq 26,  \\
Dh_Z(0)=1,  \quad Dh_Z(1)=4, \quad  Dh_Z(2)\geq 8 \\ 
Dh_Z(3) + Dh_Z(4)+ Dh_Z(5)  \geq  Dh_Z(0)+Dh_Z(1)+ Dh_Z(2)  \geq13.\\ 
\end{array} $$
The previous conditions give a contradiction if $\ell(B)<13$. They imply $Dh_Z(0)+Dh_Z(1)+ Dh_Z(2)=13$
and since $Dh_Z(0)=1$, $Dh_Z(1)= 4$, then necessarily $h_Z(2)=h_A(2)=13$. From $CB(5)$ we also get 
$Dh_Z(5)+Dh_Z(4)\geq Dh_Z(0)+Dh_Z(1)=5$. Thus $h_Z(3)\leq 21$.
\end{proof}

Putting together Proposition \ref{capnon0} and Proposition \ref{hilbdisj} we get:

\begin{thm0} \label{rango13}
Let $A\subset\Pj^4$ be a set of $13$ points satisfying conditions $(i), (ii)$, $(iii), (iv')$. Let $T$ be a form in the span
of $v_4(A)$, such that $A$ is non-redundant for $T$. Then $T$ has rank $13$.
\end{thm0}

Next, we want to find conditions for the identifiability of a quartic form in five variables, of rank 13. The following
example shows how one can construct non-identifiable forms in the span of $v_4(A)$, even if $A$ is  general.

\begin{exa0} \label{10pt} Let $A$ be a general set  of $13$ points in $\Pj^4$. \\
From the resolution of the ideal $I_A$
illustrated by sequence \eqref{eq:idD}, we know that $A$ is contained in the complete intersection of two quadric hypersurfaces. By Bertini's 
Theorem for irreducibility, the two quadrics intersect in a quartic surface $S$, whose section with a general hyperplane determined
by the linear form $\Lambda$ is an elliptic normal curve $\Gamma$ of degree $4$ in $\Pj^3$. \\
Fix a general set $W$ of $10$ points in $\Gamma$. $W$ is the residue in $\Gamma$ of two points in a complete intersection
of type $(2,2,3)$. Since $2$ points determine a linear series $g^1_2$ on $\Gamma$, then $W$ it is contained 
in a pencil of cubic surfaces of $\Pj^3$. The $h$-vector of $W$ in $\Pj^3$ is $(1,3,4,2)$. Thus $W$ is separated by cubics. It follows that
$(1,3,4,2)$ is also the  $h$-vector in $\Pj^4$ of $W$, (which is a set of points in the hyperplane $\Lambda$).\\
There  are two cubic surfaces in $\Pj^3$ which intersect $\Gamma$ properly and contain $W$. The two cubics lift to two
cubic hypersurfaces in $\Pj^4$ which contain $A$. Namely $A$ is separated by quadrics, thus
the restriction map $(I_A)_3\to (R/\Lambda)_3$ surjects
(see Proposition \ref{restri}). It follows that   $W'=A\cup W$ is contained in a complete intersection of type $(2,2,3,3)$.
The residue is a finite scheme $B$ of length $36-13-10=13$.\\
We want to prove that for a general choice of $A$, $\Lambda$, and $W$ the residue $B$ is a set of $13$ points, which determines a second
decomposition of a form in the span of $v_4(A)$.\\ 
To get the result, we claim first that $A\cup W$ is separated by cubics. Indeed in $\Pj^3$ the ideal of $W$ is generated by two quadrics and 
two cubics, because $W$ is linked $(2,2,3)$ to a pair of points. Moreover the four generators have no syzygies of degree $3$. 
The two quadrics and the two cubics spanning the ideal of $W$ in $\Pj^3$
lift to two quadrics $Q_1,Q_2$ and two cubics $K_1,K_2$, all  containing $A$, as explained above. Thus any cubic $K$ containing $A\cup W$ is of the form
$$K = \Lambda Q +Q_1\Lambda_1 + Q_2 \Lambda_2 + c_1K_1 + c_2K_2$$
where $Q$ is a quadric, $\Lambda_1,\Lambda_2$ are linear forms and $c_1,c_2\in\C$. Since $Q_1,Q_2,K_1,K_2$ contain $A$ and $\Lambda$
is generic, then $Q$ contains $A$. Then $Q$ belongs to the ideal of $A$, which is generated in degree $2$ by $Q_1,Q_2$. Thus
the ideal of $A\cup W$ is generated by  $Q_1,Q_2,K_1,K_2$ up to degree $3$. Since the four forms cannot have syzygies in degree $3$, then
one computes $h_{A\cup W}(3)=23$, and the claim holds.\\
When $A$ and $W$ are general, by semi-continuity a general set $B$ linked to $A\cup W$ by a complete intersection 
of type $2,2,3,3$ is smooth and disjoint from $A\cup W$, hence it consists of $13$ distinct points.\\
Since the $h$-vector of a complete intersection $(2,2,3,3)$ is $(1,4,8,10,8,4,1)$, as one can easily compute, then
the $h$-vector of $Z=A\cup B$, which is linked to $W$, is $(1,4,8,8,4,1)$(see formula \eqref{hlink}).
Thus  $A\cup B$ does not impose independent conditions to quartics. The intersection of the spans of $v_4(A)$ and $v_4(B)$ is a point $T$.
Experiments on numerical data prove that, for a general choice of $A$, $W$ and the linkage,  then $A$, $B$ are two non-redundant 
decompositions of $T$ (see the file {\tt ancillary.txt}). Thus $T$ has rank $13$ (by Theorem \ref{rango13})
and two different minimal decompositions. 
\end{exa0}

In the previous example, notice that $I_A$ and $I_{A\cup W}$ coincide in degree $2$, and $A\cup W$ is separated by cubics.
Thus  the Hilbert function  of $A\cup W$ and its difference are given by:
\begin{tabular}{c|cccccc}
$j$ & $0$ & $1$ & $2$ &   $3$ & $4$ &   $\dots$ \\  \hline
$h_{A\cup W}(j)$ & $1$ & $5$ &   $13$ &   $23$ &  $23$ & $\dots$ \cr
$Dh_{A \cup W}(j)$ & $1$ & $4$ &   $8$ &   $10$ & $0$ & $\dots$ 
\end{tabular}.

It follows by the mapping cone that $B$ has an $h$-vector $(1,4,8)$ identical to the $h$-vector of $A$.

\medskip

We want to prove how we can detect tensors in the span of $v_4(A)$, for $A$ general of length $\ell(A)=13$, with a second decomposition $B$
of length $13$. After Proposition \ref{hilbdisj}, we can limit ourselves to the case in which $A\cap B$ is empty. 

Call $Z$ the union $Z=A\cup B$. $Z$ has cardinality $26$ and, by Theorem \ref{CBconseq}, it satisfies $CB(4)$.

\begin{rem0} \label{Serrespec} When $A$ is general, then $I_A$ has two generators of degree $2$ and the two hypersurfaces intersect in
a smooth quartic surface $S$, whose geometry is outlined in Remark \ref{esse}. By Proposition \ref{hilbdisj}, $S$ contains $Z$. \\
Since $Z$ satisfies $CB(4)$, the Serre construction (see Section \ref{sec:Serre}) provides a rank $2$ vector bundle $\Ec$ with Chern classes
$c_1=4H - K_S=5H$ and $c_2=26$ on $S$, such that $Z$ is the zero-locus of a global section of $\Ec$. The ideal sheaf $\Ic_{Z,S}$ of $Z$ in $S$
sits in the exact sequence
\begin{equation}\label{seq}  0 \to \Oc_S \to \Ec \to \Ic_{Z,S}(5H)\to 0.
\end{equation}
Since $Z$ is contained in $14$ independent cubics, while there are only $10$ independent cubics containing $S$, then $\Ec(-2H)$ has global sections.
Indeed $h^0(\Ic_{Z,S}(4H))\geq 4$, so one has
\begin{equation} \label{h0E}
h^0\Ec(-2H)\geq 4.
\end{equation}

One computes that the Chern classes of $\Ec(-2H)$ are $c_1(\Ec(-2H))=c_1(\Ec)-4H=H$ and $c_2(\Ec(-2H))=c_2(\Ec)-2H\cdot c_1(\Ec)+2H\cdot 2H=2$.
\end{rem0}

We want to prove that $\Ec(-2H)$ has section which vanish in codimension $2$. This is true, by Proposition \ref{codim2},
unless all the global sections of  $\Ec(-2H)$ are obtained by taking one global section of a twist $h^0(\Ec(-2H-D))>0$, with $D$ non-trivial and effective,
and multiplying it by elements of $\Oc_S(D)$. We get the result by exploiting the possibilities for $D$.

\begin{rem0} Let $D=ah+\sum_{i=1}^5 b_ie_i$ be a divisor in $S$. Then  $h^0(\Oc_S(ah))\geq h^0(\Oc_S(D))$.\\
 Indeed clearly if $b_i\leq 0$ then  $h^0(\Oc_S(D-b_ie_i)\geq h^0(\Oc_S(D))$ because $-b_ie_i$ is effective. If $b_i>0$ then
 $De_i<0$, so that $e_i$ is a fixed component in $D$, hence $h^0(\Oc_S(D-e_i)= h^0(\Oc_S(D)))$ and one obtains the claim,
 arguing by induction.
\end{rem0}

\begin{prop0} \label{cod2} A general section of $\Ec(-2H)$ vanishes in a zero-dimensional scheme of length $2$.
\end{prop0} 
\begin{proof} We need to prove that there are no effective divisors $D$ such that every section of $\Ec(-2H)$ 
is a product of a fixed section $s$ of $\Ec(-2H-D)$ by divisors in $\Oc_S(D)$. By \eqref{h0E} the claim is obvious when $h^0(\Oc_S(D))<4$. \\
Write $D=ah+\sum_{i=1}^5 b_ie_i$ and order the $b_i$'s so that $b_1\leq\dots\leq b_5$.\\
Since $h^0(\Oc_S(h))=3$, by the previous remark we can exclude the cases in which $a\leq 1$.\\
Assume $a=2$. Since $h^0(\Oc_S(2h))=6$, we have $b_1,b_2\geq -1$, $b_3=b_4=b_5\geq 0$. From $h^0(\Ec(-2H-D))>0$ and from\eqref{seq} we get
that $Z$ belongs to a divisor of type $3H-D=7h-\sum_{i=1}^5 (b_i+3)e_i$. From the previous inequalities we get
$h^0(\Oc_S(3H-D))\leq h^0(\Oc_S(7h-2e_1-2e_2-3e_3-3e_4-3e_5))=36-3-3-6-6-6=12$. It follows that divisors of type $3H-D$
cannot contain a general set $A$ of $13$ points of $S$, thus they cannot contain $Z$, a contradiction.\\
Assume $a=3$. Since $h^0(\Oc_S(3h))=10$, then $(b_1,\dots,b_5)$ is greater or equal than the following $5$-tuples:
\begin{itemize}
\item[1)] $(-3,0,0,0,0)$,
\item[2)] $(-2,-2,0,0,0)$,
\item[3)] $(-2,-1,-1,-1,0)$,
\item[4)] $(-1,-1,-1,-1,-1).$
\end{itemize}
In case 1) one computes that  $h^0(\Oc_S(3H-D))\leq h^0(\Oc_S(6h-3e_2-3e_3-3e_4-3e_5))=28-6-6-6-6=4$. In case 2) one computes that 
$h^0(\Oc_S(3H-D)) \leq h^0(\Oc_S(6h-e_1-e_2-3e_3-3e_4-3e_5))=28-1-1-6-6-6=8$. In case 3) one computes that 
$h^0(\Oc_S(3H-D)) \leq h^0(\Oc_S(6h-e_1-2e_2-2e_3-2e_4-3e_5))=28-1-3-3-3-6=12$.  In case 4) one computes that 
$h^0(\Oc_S(3H-D)) \leq h^0(\Oc_S(6h-2e_1-2e_2-2e_3-2e_4-2e_5))=28-3-3-3-3-3=13$. In any event $h^0(\Oc_S(3H-D))\leq 13$.
Thus divisors of type $3H-D$
cannot contain a general set $A$ of $13$ points of $S$, thus they cannot contain $Z$, a contradiction.\\
Assume $a=4$. Since the general set $A$ of $13$ points is contained in a divisor of type $3H-D=5h-\sum_{i=1}^5(b_i+3)e_i$,
so that $h^0(\Oc_S(3H-D))\geq 14$, and since $h^0(\Oc_S(5h))=21$, then $(b_1+3,\dots,b_5+3)$  is greater or equal than the following $5$-tuples:
\begin{itemize}
\item[1)] $(-3,-1,0,0,0)$,
\item[2)] $(-2,-2,-1,0,0)$,
\item[3)] $(-2,-1,-1,-1,-1)$.
\end{itemize}
In any case, one computes that  $D$ cannot pass through $13$ general points, a contradiction. \\
The cases in which $a\geq 5$ are similar but much easier.\\
It follows that a general section of $\Ec(-2H)$ vanishes in a finite set $Z'$ of length 
$$\ell(Z')=c_2(\Ec(-2H))=c_2(\Ec) -2H\cdot c_1(\Ec) + (2H)^2 = 26-10H^2+4H^2=26-24=2.$$
\end{proof}

It follows that $Z$ and $Z'$ are zero-loci of sections of two different twists of $\Ec$.
The vanishing loci of two different twists of $\Ec$ are connected in a double linkage. We explain the procedure in one example.

\begin{exa0} \label{mainex} Consider a general set $A$ of $13$ points in $\Pj^4$ and take a set $Z'$ of $2$ general points in the quartic surface $S$,
complete intersection of type $1,1,1,2$, determined by $A$. The ideal of $Z'$ has the Koszul resolution:
\begin{multline*}
0 \rightarrow R(-5) \rightarrow R(-3)\oplus R(-4)^{\oplus 3} \rightarrow \oplus R(-3)^{\oplus 3}\oplus R(-2)^{\oplus 3}\rightarrow\\
\rightarrow R(-2) \oplus R(-1)^{\oplus 3} \rightarrow I_{Z'} \rightarrow 0.
\end{multline*}  
Fix a general hyperplane $\pi$ through $Z'$ and consider a general cubic $F$ which contains $Z'$ and $A$. Since $Z',A$ are general, for a general choice
of $F$ the curve $F\cap S$ is irreducible. $F$ exists since $h_A(3)=13$ so $(I_A)_3$ has
dimension $ 7$. Link $Z'$ to a set $Z''$ by the intersection of $S,\pi,F$. Since $A$ is general then $F$ is a general cubic through $Z'$ and the two quadrics which
determine $S$ are general quadrics through $Z'$. Being $I_{Z'}$ generated in degree $2$,  we get that $Z''$ is a set of $10$ points in the hyperplane $\pi$.
A resolution of $I_{Z''}$ is determined by the mapping cone associated to the diagram
{\small$$
\begin{matrix}
0 &\to& R(-8) &\to& \begin{matrix}R(-7) \\ \oplus \\ R(-6)^{\oplus 2} \\ \oplus \\ R(-5) \end{matrix}&\to& 
\begin{matrix}R(-5)^{\oplus 2}\\ \oplus \\ R(-4)^{\oplus 2}\\ \oplus \\ R(-3)^{\oplus 2} \end{matrix} &\to& 
\begin{matrix} R(-3) \\ \oplus \\ R(-2)^{\oplus 2}\\ \oplus \\ R(-1)\end{matrix}  &\to& I_{W'} &\to& 0 \\ 
& & \vspace{.1cm} & & & & & & & &  & &  \\

 & & \big\downarrow & & \big\downarrow & &  \big\downarrow & & \big\downarrow & & \big\downarrow  & &  \\

& & \vspace{.1cm} & & & & & & & &  & &  \\

0 &\to& R(-5) &\to& \begin{matrix}R(-4)^{\oplus 3} \\ \oplus \\ R(-3) \end{matrix}&\to& 
\begin{matrix}R(-2)^{\oplus 3}\\ \oplus \\ R(-3)^{\oplus 3}\end{matrix} &\to& \begin{matrix}R(-2)\\ \oplus \\ R(-1)^{\oplus 3}\end{matrix}  &\to& I_{Z'} &\to& 0
\end{matrix}$$}
where $W'=S\cap F\cap \pi$. We get
\begin{multline*}
0 \to R(-7)^{\oplus 2} \rightarrow R(-6)^{\oplus 4}\oplus R(-5)^{\oplus 3} \to R(-5)^{\oplus 2}\oplus R(-4)^{\oplus 5}\oplus R(-3)^{\oplus 2}\to\\
\to  R(-3)^{\oplus 2}\oplus R(-2)^{\oplus 2} \oplus R(-1) \to I_{Z''} \to 0.
\end{multline*}
Notice from the resolution that $(I_{Z''})_3$ has dimension $27$. Since the linear space of cubics containing $S$ is $10$-dimensional, then
there exists a space of dimension at least $27-10-13= 4$ of cubics containing $A$ and not containing $S$. These cubics determine a linear system
of divisors $3H$ on $S$ which contain $Z''\cup A$. Since $F\cap S$ is irreducible, for a general choice of a cubic $F'$ through $Z''\cup A$
the intersection $S\cap F\cap F''$ is a finite set $W''$ of length $4\cdot 3\cdot 3=36$.
Linking $Z''$ with $S,F,F''$ we obtain thus a set $Z$ of length $26$ containing $A$. Since $A, Z''$ are general and the ideal of $Z'$ is generated in degree $2$, 
then by \cite{PeskineSzpiro74} we may assume that $Z$ is a set of $26$ distinct points. The resolution of $Z$ comes from the mapping cone of the diagram 
{\small $$
\begin{matrix}
0 &\to& R(-10) &\to& \begin{matrix} R(-8)^{\oplus 2} \\ \oplus \\ R(-7)^{\oplus 2}  \end{matrix}&\to& 
\begin{matrix}R(-6)\\ \oplus \\ R(-5)^{\oplus 4}\\ \oplus \\ R(-4) \end{matrix} &\to& 
\begin{matrix} R(-3)^{\oplus 2} \\ \oplus \\ R(-2)^{\oplus 2} \end{matrix}  &\to& I_{W''} &\to& 0 \\ 
& & \vspace{.1cm} & & & & & & & &  & &  \\

 & & \big\downarrow & & \big\downarrow & &  \big\downarrow & & \big\downarrow & & \big\downarrow  & &  \\

& & \vspace{.1cm} & & & & & & & &  & &  \\

0 &\to& R(-7)^{\oplus 2} &\to& \begin{matrix}R(-6)^{\oplus 4} \\ \oplus \\ R(-5)^{\oplus 3}  \end{matrix}&\to& 
\begin{matrix}R(-5)^{\oplus 2}\\ \oplus \\ R(-4)^{\oplus 5}\\ \oplus \\ R(-3)^{\oplus 2}\end{matrix} &\to& \begin{matrix}R(-3)^{\oplus 2}\\ \oplus \\ R(-2)^{\oplus 2}\\ \oplus \\ R(-1)\end{matrix}  &\to& I_{Z''} &\to& 0
\end{matrix}$$}
so that we obtain
\begin{multline*}
0 \to R(-9) \rightarrow R(-7)^{\oplus 2} R(-6)^{\oplus 4}\to  R(-5)^{\oplus 5}\oplus R(-4)^{\oplus 5}\to\\
\to  R(-3)^{\oplus 4}\oplus R(-2)^{\oplus 2}  \to I_{Z} \to 0.
\end{multline*}
Thus the $h$-vector of $Z$ is $(1,4,8,8,4,1)$. Since the last module of the resolution is $1$-dimensional, the
set $Z$ is arithmetically Gorenstein, by Proposition \ref{DGO}. In particular, the resolution is auto-dual.
The residue of $A$ in $Z$ is a set $B$ of $13$ points. \\
Since $h_Z(4)=25=26-1$, the spaces
$\langle v_4(A)\rangle$ and $\langle v_4(B) \rangle$ meet in exactly one point $T$, which represents a form with $2$
different decompositions of length $13$. 
\end{exa0}

\begin{rem0} In the previous example, $A$ and $B$ are linked by the arithmetically Gorenstein scheme $Z$.
Thus, by the Gorenstein liaison procedure (see Section \ref{sec:liaison}), a resolution of the ideal of $B$ comes out from the resolution of the ideals $I_Z$ and $I_A$,
by taking the mapping cone of the diagram
{\small $$
\begin{matrix}
0 &\to& R(-9) &\to& \begin{matrix} R(-7)^{\oplus 2} \\ \oplus \\ R(-6)^{\oplus 4}  \end{matrix}&\to& 
\begin{matrix}R(-5)^{\oplus 5}\\ \oplus \\ R(-4)^{\oplus 5} \end{matrix} &\to& 
\begin{matrix} R(-3)^{\oplus 4} \\ \oplus \\ R(-2)^{\oplus 2} \end{matrix}  &\to& I_Z &\to& 0 \\ 
& & \vspace{.1cm} & & & & & & & &  & &  \\

 & & \big\downarrow & & \big\downarrow & &  \big\downarrow & & \big\downarrow & & \big\downarrow  & &  \\

& & \vspace{.1cm} & & & & & & & &  & &  \\

0 &\to& R(-6)^{\oplus 8} &\to& R(-5)^{\oplus 28} &\to& R(-4)^{\oplus 33} &\to& \begin{matrix}R(-3)^{\oplus 12}\\ \oplus \\ R(-2)^{\oplus 2}\end{matrix}  &\to& I_A &\to& 0
\end{matrix}$$}
One realizes immediately that the Betti numbers of a resolution of $B$ coincide with the Betti numbers of a resolution of $A$.
\end{rem0}

The following proposition is an exercise of duality.

\begin{prop0} Assume $A$ is a general set of $13$ points in $\Pj^4$ and let $T$ be a form in the span of $v_4(A)$ such that $A$ is a 
non-redundant decomposition of  $T$. If $T$ has a a second decomposition $B$ of length $13$ with $A\cap B=\emptyset$, then $B$ is obtained from $A$
by a double linkage as illustrated in Example  \ref{mainex}, starting with a non-necessarily reduced scheme $Z'$ of length $2$ in $S$. In particular,
the set $Z=A\cup B$ has $h$-vector $(1,4,8,8,4,1)$.
\end{prop0}
\begin{proof} We know that two general quadrics containing $A$ intersect in a smooth irreducible quartic surface $S$. Let $B$ be another minimal
decomposition of $T$, disjoint from $A$, and write $Z=A\cup B$. By Proposition \ref{hilbdisj} $Z$ is contained in $S$. By Theorem \ref{CBconseq}, $Z$ satisfies $CB(4)$,
hence it is associated to a rank $2$ vector bundle $\Ec$ on $S$. By Proposition \ref{cod2} $\Ec(-2H)$ has a section vanishing in a zero-dimensional scheme $Z'$
of length $2$. It is classically known that schemes of length $2$ on a smooth surface $S$ form an irreducible family whose general element is a set of $2$
distinct points. Any scheme of length $2$ is separated by hyperplanes, thus it is aligned and complete intersection of type $1,1,1,2$.
If $Z'$ is not smooth, then it consists of a length $2$ scheme concentrated at a point $P\in S$ and in a tangent line $L$ to $S$ at $P$. Any tangent
line $L$ uniquely determines  $Z'$.\\
Let $s$ be a global section of $\Ec$ which vanishes at $Z$ and let $s'$ be a section of $\Ec(-2H)$ vanishing at $Z'$. By the exact sequence
$$ 0\to \Oc_S \to \Ec(-2H) \to \Ic_{Z'}(H) \to 0,$$
where $\Ic_{Z'}$ is the ideal sheaf of $Z'$ on $S$, the section $s$ determines a section of $\Ic_{Z'}(3H)$, i.e. a divisor $D_F$ of type $3H$
which vanishes where $s\wedge s'$ vanishes. Hence $D_F$ lifts to a cubic $F$ containing $Z'\cup Z$. Fix a general hyperplane $\pi$ containing $Z'$, 
which determines a section of $\Ec(-2H)$. Call $Z''$ the residual scheme $Z''$ of $Z'$ with respect to the intersection $ F\cap \pi$ on $S$.
The mapping cone of the diagram
$$
\begin{matrix}
0 &\to& \Oc_S(-4H) &\to& \Oc_S(-3H)\oplus \Oc_S(-H)   &\to&\Ic_{F\cap \pi\cap S} &\to& 0 \\ 
 & & \big\downarrow & & \big\downarrow & &  \big\downarrow & &   \\
0 &\to& \Oc_S(-H) &\to& \Ec(-3H)  &\to& \Ic_{Z'} &\to& 0
\end{matrix}$$
determines the following exact sequence for the ideal sheaf $\Ic_{Z''}$ of $Z''$ on $S$
\begin{equation} \label{zz}
 0\to \Ec^\vee (-H) \to \Oc_S(-3H) \oplus \Oc_S(-3H)\oplus \Oc_S(-H) \to \Ic_{Z''} \to 0.
 \end{equation}
Notice that, in the diagram, the map $\Ec^\vee (-H) \to \Oc_S(-3H)$ is the dual of the map $ \Oc_S(-H)\to \Ec(-H)$ defined by $s$.
Furthermore, the  central vertical map corresponds to the section $s$ of $\Ec$ and a section $s_0$ of $\Ec(-2H)$,
where $s_0\wedge s'$ defines the divisor $\pi\cap S$. The map $\Ec^\vee (-H) \to \Oc_S(-3H) \oplus \Oc_S(-3H)\oplus \Oc_S(-H) $
in the mapping cone is  the dual of the map $\Oc_S(-3H) \oplus \Oc_S(-3H)\oplus \Oc_S(-H)\to \Ec$ defined by $s',s_0,s$.
Thus the three divisors, of type $3H,3H,H$ containing  $Z''$, defined by the map  $\Oc_S(-3H) \oplus \Oc_S(-3H)\oplus \Oc_S(-H) \to \Ic_{Z''}$,
correspond to the zero-loci of the wedge products $s'\wedge s, s_0\wedge s, s'\wedge s_0$ respectively.
Let $F'$ be a cubic in $\Pj^4$ such that $F\cap S$ corresponds to the zero-locus of $s_0\wedge s$. Clearly $F'$ contains both
$Z$ and $Z''$. Sequence \eqref{zz} proves that $Z''$ is cut set theoretically by $F,F',\pi,S$, thus $F\cap F'\cap \pi$ is zero-dimensional.
The mapping cone of the diagram 
$$
\begin{matrix}
0 &\to& \Oc_S(-6H) &\to& \Oc_S(-3H)\oplus \Oc_S(-3H)   &\to&\Ic_{F\cap F'\cap S} &\to& 0 \\ 
 & & \big\downarrow & & \big\downarrow & &  \big\downarrow & &   \\
0 &\to& \Ec^\vee(-H) &\to& \Oc_S(-3H) \oplus \Oc_S(-3H)\oplus \Oc_S(-H) &\to& \Ic_{Z''} &\to& 0
\end{matrix}$$
shows that the residue of $Z''$ with respect to $F\cap F'\cap S$ is exactly the zero-locus of the section defined 
by the dual of the map $\Ec^\vee(-H) \to \Oc_S(-H)$. Hence the residue is $Z$.\\
When $Z'$ is a set of two distinct points, we get back exactly the double linkage described in Example \ref{mainex}.
Since every set $Z'$ of length $2$ is a limit of a family of reduced sets $\Zc'$ and every hyperplane containing $Z'$ is a limit
of a family of hyperplanes containing $\Zc'$, then, even if $Z'$ is non-reduced, the scheme $Z$ is a limit of schemes
obtained as in Example  \ref{mainex}.\\
For any scheme $Z'$ of length $2$ in $S$ we have $h^0(\Ic_{Z'}(2H))=11$. Thus $h^0(\Ec(-H))=16$, so that
$h^0(\Ic_Z(4H))=16$. Hence $\dim(I_Z)_4=h^0(\Ic_Z(4H))+\dim(I_S)_4=16+29=45$. It follows that
$h_Z(4)=70-45=25=\ell(Z)-1$. Hence $ \langle v_4(A)\rangle\cap\langle v_4(B)\rangle = \{T\} $, i.e. $T$ is uniquely determined by $B$.\\
The claim follows.
\end{proof}

Next results explain how one can parameterize sets $B$ such that the spans of $v_4(A),v_4(B)$ meet at a point. 

\begin{prop0} Fix a general set $A$ of $13$ points in $\Pj^4$ and let $S$ be the (smooth) quartic surface,
intersection of the quadrics containing $A$. Let $C$ be a general hyperplane section of $S$. Then for a general choice 
of a set $Z''$ of $10$ points on $C$ the $h$-vector of $A\cup Z''$ is $(1,4,8,10)$. Thus, the residue $B$ of $A\cup Z''$
is a set such that the spans of $v_4(A),v_4(B)$ meet at a point $T$.
\end{prop0}
\begin{proof} Adding to $A$ $10$ points in some hyperplane section of $S$ obtained as in Example \ref{mainex}, we get a set of points whose
residue with respect to a complete intersection of type $2,2,3,3$ is a set $B$ as in the previous proposition. Since the residue of such $B$
has $h$-vector $(1,4,8,10)$, the same holds for the addition to $A$ of $10$ general points in $C$. The second
claim follows immediately by the mapping cone.
\end{proof}

It remains to prove that for a general choice of a set $Z''$ of $10$ points in a hyperplane section of $S$ the residue $B$ of $A\cup Z''$ in a complete
intersection of $S$ and two cubics defines a quartic $T=\langle v_4(A)\rangle\cap\langle v_4(B)\rangle$ such that $A$ is
a non-redundant decomposition for $T$. \\
In practice, this amounts to find the expression of a general such $T$ in terms of coordinates of the points of $A$ and
show that no coefficients of the expression vanish.\\
It is clear that the previous property is open, thus it is sufficient to check it in one example.

\begin{exa0} \label{nonredund} Fix a general set $A$ of $13$ points in $\Pj^4$,  which determines a general (smooth) surface $S$,
complete intersection of two quadrics. Fix a general set $Z''$ of $10$ points on a general hyperplane section of $S$.
From the coordinates of the set $A$ one computes the ideal $I_A$. From the coordinates of the set $A\cup Z''$ and the mapping cone
procedure, one computes the ideal $I_B$. If everything is general enough, then $B\cap A$ is empty and $h_Z(4)=25$, where $Z=A\cup B$. 
It follows that in the polynomial ring $R=\C [x_0,\dots,x_4]$ the subspaces $(I_A)_4$ and $(I_B)_4$ of $R_4$ are both $57$-dimensional,
but their sum $U=(I_A)_4+(I_B)_4$ is a linear subspace of codimension $1$ in $R_4$. Coefficients for the equation of $U$
provide coefficients for the form $T$ in the intersection of $\langle v_4(A)\rangle$ and $\langle v_4(B)\rangle$. We refer to \cite{AngeC20} for 
details.\\
Thus, from the coordinates of $A$ and $Z''$ one determines $T$ explicitly, and also an expression of $T$ as combination
of rank $1$ tensors representing the points of $A$. Experiments prove the existence of $T$ for which no coefficients
vanish. We refer for that to the  file ancillary.txt.\\
If $A$ is general enough, then the span of tangent spaces to the Veronese variety $v_4(\Pj^4)$ at the points of $v_4(A)$ determines 
a linear subspace of the expected projective dimension $64$ in the projective space spanned by $v_4(\Pj^4)$. This is easy to check
directly for the points $A$ above (see the file ancillary.txt). From Remark \ref{stacked} it follows  that there are no positive-dimensional  families 
of decompositions for the form $T$.
\end{exa0} 

\begin{coro0} \label{cor10} Fix a general set $A$ of $13$ points in $\Pj^4$. Then the span of $v_4(A)$, which is a $\Pj^{12}$, contains an
irreducible  $10$-dimensional family $\Tc$ of points $T$ which have a second decomposition $B$ obtained as in Example \ref{mainex}. 
\end{coro0}
\begin{proof} The family $\Zc''$ of sets $Z''$ of $10$ points  contained in some hyperplane section of $S$ is irreducible and $14$-dimensional
and maps to the family of points $T\in\langle v_4(A)\rangle$ obtained as in Example \ref{mainex}. It is thus sufficient to prove that
the general fiber is $4$-dimensional.\\
 Indeed, by the previous example, the family of points $T$ is a finite image of the family $\Bc$ of sets $B$ linked to $A\cup Z''$, when $Z''$ moves.
For a general choice of $B$ in the family $\Bc$, the sets $Z''$ determining $B$ are obtained as the residue in a complete
intersection of $S$ and two more cubics. Thus the fiber over $B$ of the map $\Zc''\to \Bc$  is determined by the choice of a pencil  of cubics through
$A\cup B$. Since the vector space of cubics through $A\cup B$ is $4$-dimensional, the fiber corresponds to
the choice of a line in a $\Pj^3$. The claim follows.
\end{proof}

\begin{rem0} In principle, as in the case $r=12$, also for an expression $A$ of length $13$ one could try to construct an algorithm
that checks if $A$ is minimal and unique.\\
In particular, the algorithm should test if $A$ satisfies $(i), (ii), (iii)$, if for any choice of two points $P_i,P_j\in A$ the set $A\setminus\{P_i,P_j\}$
satisfies $(iv)$, and if for any choice of the index $i=1,\dots,13$ the algorithm for $r=12$ gives the identifiability of $T-T_i$.\\
Then, in order to detect if $T$ belongs to the bad family $\Tc$ or not, the algorithm should construct a parametrization of the bad family.
This step is much more delicate than the corresponding step $7$ for the case $r=12$, because it requires to parametrize
the choice of a set of $10$ points in an elliptic curve.\\
We intend to devote a future paper to practical solutions of this last technical step.
\end{rem0}

\bibliographystyle{amsplain}
\bibliography{biblioLuca2}

\providecommand{\bysame}{\leavevmode\hbox to3em{\hrulefill}\thinspace}
\providecommand{\MR}{\relax\ifhmode\unskip\space\fi MR }
\providecommand{\MRhref}[2]{%
  \href{http://www.ams.org/mathscinet-getitem?mr=#1}{#2}
}
\providecommand{\href}[2]{#2}
\begin{thebibliography}{10}

\bibitem{AlexHir95}
J.~Alexander and A.~Hirschowitz, \emph{Polynomial interpolation in several
  variables}, J. Algebraic Geom. \textbf{4} (1995), 201--222.

\bibitem{AngeBocciC18}
E.~Angelini, C.~Bocci, and L.~Chiantini, \emph{Real identifiability vs complex
  identifiability}, Lin. Multilin. Algebra \textbf{66} (2018), 1257--1267.

\bibitem{AngeC20}
E.~Angelini and L.~Chiantini, \emph{On the identifiability of ternary forms},
  Lin. Alg. Applic. \textbf{599} (2020), 36--65.

\bibitem{AngeC}
\bysame, \emph{Minimality and uniqueness for decompositions of specific ternary
  forms}, Math. of Comput. \textbf{91} (2022), 973--1006.

\bibitem{AngeCMazzon19}
E.~Angelini, L.~Chiantini, and A.~Mazzon, \emph{Identifiability for a class of
  symmetric tensors}, Mediterr. J. Math. \textbf{16} (2019), 97.

\bibitem{AngeCVan18}
E.~Angelini, L.~Chiantini, and N.~Vannieuwenhoven, \emph{Identifiability beyond
  {K}ruskal's bound for symmetric tensors of degree 4}, Rend. Lincei Mat.
  Applic. \textbf{29} (2018), 465--485.

\bibitem{Ball19}
E.~Ballico, \emph{An effective criterion for the additive decompositions of
  forms}, Rend. Ist. Matem. Trieste \textbf{51} (2019), 1--12.

\bibitem{BallBern12a}
E.~Ballico and A.~Bernardi, \emph{Decomposition of homogeneous polynomials with
  low rank}, Math. Zeit. \textbf{271} (2012), 1141--1149.

\bibitem{Brun80}
J.~Brun, \emph{Les fibres de rang deux sur {$\mathbb P^2$} et leur sections},
  Bull. Soc. Math. France \textbf{108} (1980), 457--473.

\bibitem{CCi02a}
L.~Chiantini and C.~Ciliberto, \emph{Weakly defective varieties}, Trans. Amer.
  Math. Soc. \textbf{354} (2002), 151--178.

\bibitem{COtt12}
L.~Chiantini and G.~Ottaviani, \emph{On generic identifiability of 3-tensors of
  small rank}, SIAM J. Matrix Anal. Appl. \textbf{33} (2012), 1018--1037.

\bibitem{COttVan14}
L.~Chiantini, G.~Ottaviani, and N.~Vannieuwenhoven, \emph{An algorithm for
  generic and low-rank specific identifiability of complex tensors}, SIAM J.
  Matrix Anal. Appl. \textbf{35} (2014), 1265--1287.

\bibitem{COttVan17b}
\bysame, \emph{Effective criteria for specific identifiability of tensors and
  forms}, SIAM J. Matrix Anal. Appl. \textbf{38} (2017), 656--681.

\bibitem{COttVan17a}
\bysame, \emph{On generic identifiability of symmetric tensors of subgeneric
  rank}, Trans. Amer. Math. Soc. \textbf{369} (2017), 4021--4042.

\bibitem{DavisGerOre85}
E.~Davis, A.~Geramita, and F.~Orecchia, \emph{{G}orenstein algebras and the
  {C}ayley-{B}acharach theorem}, Proc. Amer. Math. Soc. \textbf{93} (1985),
  593--597.

\bibitem{Derksen13}
H.~Dersken, \emph{Kruskal's uniqueness inequality is sharp}, Linear Alg.
  Applic. \textbf{438} (2013), 708--712.

\bibitem{Domanov}
I.~Domanov, \emph{Study of canonical polyadic decomposition of higher-order
  tensors}, Ph.D. thesis, Arenberg Doctoral School, KU Leuven, 2013.

\bibitem{DomaLath15}
I.~Domanov and L.~{De Lathauwer}, \emph{Generic uniqueness conditions for the
  canonical polyadic decomposition and {INDSCAL}}, SIAM J. Matrix Anal. Appl.
  \textbf{36} (2015), 1567--1589.

\bibitem{Eisenbud}
D.~Eisenbud, \emph{Commutative algebra, with a view towards algebraic
  geometry}, Graduate Texts in Mathematics, Springer, Berlin, New York NY,
  1995.

\bibitem{Ferrand75}
D.~Ferrand, \emph{Courbes gauches et fibres de rang {2}}, C. R. Acad. Sci.
  Paris \textbf{281} (1975), 345--347.

\bibitem{Macaulay2}
D.~Grayson and M.~Stillman, \emph{Macaulay 2, a software system for research in
  algebraic geometry}, Available online http://www.math.uiuc.edu/Macaulay2.

\bibitem{Hart78}
R.~Hartshorne, \emph{Stable vector bundles of rank {2} in {$\mathbb P^3$}},
  Math. Ann. \textbf{238} (1978), 229--280.

\bibitem{Hartshorne}
\bysame, \emph{Algebraic geometry}, Graduate Texts in Math., Springer, Berlin,
  New York NY, 1992.

\bibitem{Kruskal77}
J.B. Kruskal, \emph{Three-way arrays: rank and uniqueness of trilinear
  decompositions, with application to arithmetic complexity and statistics},
  Linear Algebra Appl. \textbf{18} (1977), 95--138.

\bibitem{Mazzon20}
A.~Mazzon, \emph{On a geometric method for the identifiability of forms},
  Boll.UMI \textbf{13} (2020), 137--154.

\bibitem{Migliore}
J.~Migliore, \emph{Introduction to liaison theory and deficiency modules},
  Progress in Mathematics, vol. 165, Birk{\"a}user, Basel, Boston MA, 1998.

\bibitem{PeskineSzpiro74}
C.~Peskine and L.~Szpiro, \emph{Liaison des vari{\'e}t{\'e}s alg{\'e}briques},
  Invent. Math. \textbf{26} (1974), 271--302.

\bibitem{RaoLiZhang18}
W.~Rao, D.~Li, and J.Q. Zhang, \emph{A tensor-based approach to {L}-shaped
  arrays processing with enhanced degrees of freedom}, IEEE Signal Proc. Lett.
  \textbf{25} (2018), 1--5.

\bibitem{Walter95}
C.~Walter, \emph{The minimal free resolution of the homogeneous ideal of s
  general points in {$\mathbb P^4$}}, Math. Zeit. \textbf{219} (1995),
  231--234.

\end{thebibliography}

\end{document}